\documentclass[12pt]{article}

\usepackage{amsfonts,amssymb,amsmath,subeqnarray,eqsection,indent}
\usepackage{bm}    
\usepackage{cite}  

\date{February 24, 2008 \\[1mm] revised February 4, 2009 \\[1cm]
      {\em Dedicated to the memory of Pierre Leroux}}

\begin{document}

\title{\vspace*{-2cm}
       Some variants of the exponential formula, \\
       with application to the \\
       multivariate Tutte polynomial
       (alias Potts model)}

\author{
     \\[-5mm]
     {\small Alexander D.~Scott} \\[-2mm]
     {\small\it Mathematical Institute}  \\[-2mm]
     {\small\it University of Oxford} \\[-2mm]
     {\small\it 24--29 St.~Giles}\/ \\[-2mm]
     {\small\it Oxford OX1 3LB, England}                         \\[-2mm]
     {\small\tt scott@maths.ox.ac.uk}                        \\[5mm]
     {\small Alan D.~Sokal\thanks{Also at Department of Mathematics,
           University College London, London WC1E 6BT, England.}}  \\[-2mm]
     {\small\it Department of Physics}       \\[-2mm]
     {\small\it New York University}         \\[-2mm]
     {\small\it 4 Washington Place}          \\[-2mm]
     {\small\it New York, NY 10003 USA}      \\[-2mm]
     {\small\tt sokal@nyu.edu}               \\[-2mm]
     {\protect\makebox[5in]{\quad}}  
     \\
}

\maketitle
\thispagestyle{empty}   

\begin{abstract}
We prove some variants of the exponential formula
and apply them to the multivariate Tutte polynomials
(also known as Potts-model partition functions) of graphs.
We also prove some further identities for the multivariate Tutte polynomial,
which generalize an identity for counting connected graphs
found by Riordan, Nijenhuis, Wilf and Kreweras
and in more general form by Leroux and Gessel,
and an identity for the inversion enumerator of trees
found by Mallows, Riordan and Kreweras.
Finally, we prove a generalization of M\"obius inversion
on the partition lattice.
\end{abstract}


\bigskip
\noindent
{\bf Mathematics Subject Classification (MSC 2000) codes:}
05A19 (Primary);
05A15, 05A18, 05C15, 05C30, 05C99, 06A07, 82B20 (Secondary).

\clearpage

\newtheorem{theorem}{Theorem}[section]
\newtheorem{proposition}[theorem]{Proposition}
\newtheorem{lemma}[theorem]{Lemma}
\newtheorem{corollary}[theorem]{Corollary}
\newtheorem{definition}[theorem]{Definition}
\newtheorem{conjecture}[theorem]{Conjecture}
\newtheorem{question}[theorem]{Question}
\newtheorem{example}[theorem]{Example}

\renewcommand{\theenumi}{\alph{enumi}}
\renewcommand{\labelenumi}{(\theenumi)}
\def\eop{\hbox{\kern1pt\vrule height6pt width4pt
depth1pt\kern1pt}\medskip}
\def\prf{\par\noindent{\bf Proof.\enspace}\rm}
\def\rmk{\par\medskip\noindent{\bf Remark\enspace}\rm}

\newcommand{\be}{\begin{equation}}
\newcommand{\ee}{\end{equation}}
\newcommand{\<}{\langle}
\renewcommand{\>}{\rangle}
\newcommand{\widebar}{\overline}
\def\reff#1{(\protect\ref{#1})}
\def\spose#1{\hbox to 0pt{#1\hss}}
\def\ltapprox{\mathrel{\spose{\lower 3pt\hbox{$\mathchar"218$}}
    \raise 2.0pt\hbox{$\mathchar"13C$}}}
\def\gtapprox{\mathrel{\spose{\lower 3pt\hbox{$\mathchar"218$}}
    \raise 2.0pt\hbox{$\mathchar"13E$}}}
\def\textprime{${}^\prime$}
\def\proof{\par\medskip\noindent{\sc Proof.\ }}
\def\firstproof{\par\medskip\noindent{\sc First Proof.\ }}
\def\secondproof{\par\medskip\noindent{\sc Second Proof.\ }}
\def\thirdproof{\par\medskip\noindent{\sc Third Proof.\ }}
\def\qed{ $\square$ \bigskip}
\def\proofof#1{\bigskip\noindent{\sc Proof of #1.\ }}
\def\firstproofof#1{\bigskip\noindent{\sc First Proof of #1.\ }}
\def\secondproofof#1{\bigskip\noindent{\sc Second Proof of #1.\ }}
\def\thirdproofof#1{\bigskip\noindent{\sc Third Proof of #1.\ }}
\def\altproofof#1{\bigskip\noindent{\sc Alternate Proof of #1.\ }}
\def\half{ {1 \over 2} }
\def\third{ {1 \over 3} }
\def\twothird{ {2 \over 3} }
\def\smfrac#1#2{\textstyle{#1\over #2}}
\def\smhalf{ \smfrac{1}{2} }
\newcommand{\real}{\mathop{\rm Re}\nolimits}
\renewcommand{\Re}{\mathop{\rm Re}\nolimits}
\newcommand{\imag}{\mathop{\rm Im}\nolimits}
\renewcommand{\Im}{\mathop{\rm Im}\nolimits}
\newcommand{\sgn}{\mathop{\rm sgn}\nolimits}
\newcommand{\tr}{\mathop{\rm tr}\nolimits}
\newcommand{\supp}{\mathop{\rm supp}\nolimits}
\def\hboxscript#1{ {\hbox{\scriptsize\em #1}} }
\renewcommand{\emptyset}{\varnothing}

\newcommand{\restrict}{\upharpoonright}
\renewcommand{\implies}{\;\Longrightarrow\;}

\newcommand{\scra}{{\mathcal{A}}}
\newcommand{\scrb}{{\mathcal{B}}}
\newcommand{\scrc}{{\mathcal{C}}}
\newcommand{\scrf}{{\mathcal{F}}}
\newcommand{\scrg}{{\mathcal{G}}}
\newcommand{\scrh}{{\mathcal{H}}}
\newcommand{\scrk}{{\mathcal{K}}}
\newcommand{\scrl}{{\mathcal{L}}}
\newcommand{\scro}{{\mathcal{O}}}
\newcommand{\scrp}{{\mathcal{P}}}
\newcommand{\scrr}{{\mathcal{R}}}
\newcommand{\scrs}{{\mathcal{S}}}
\newcommand{\scrt}{{\mathcal{T}}}
\newcommand{\scrv}{{\mathcal{V}}}
\newcommand{\scrw}{{\mathcal{W}}}
\newcommand{\scrz}{{\mathcal{Z}}}

\newcommand{\ahat}{{\widehat{a}}}
\newcommand{\Zhat}{{\widehat{Z}}}
\renewcommand{\k}{{\mathbf{k}}}
\newcommand{\n}{{\mathbf{n}}}
\newcommand{\ttt}{{\mathbf{t}}}
\newcommand{\vv}{{\mathbf{v}}}
\newcommand{\bv}{{\mathbf{v}}}
\newcommand{\w}{{\mathbf{w}}}
\newcommand{\x}{{\mathbf{x}}}
\newcommand{\cc}{{\mathbf{c}}}
\newcommand{\zero}{{\mathbf{0}}}
\newcommand{\one}{{\mathbf{1}}}
\newcommand{\bdelta}{{\boldsymbol{\delta}}}
\newcommand{\bpi}{{\bm{\pi}}}
\newcommand{\B}{{\rm\bm{B}}}

\newcommand{\C}{{\mathbb C}}
\newcommand{\Z}{{\mathbb Z}}
\newcommand{\N}{{\mathbb N}}
\newcommand{\Q}{{\mathbb Q}}
\newcommand{\R}{{\mathbb R}}
\newcommand{\RR}{{\mathbb R}}

%
%
\newcommand{\stirlingsubset}[2]{\genfrac{\{}{\}}{0pt}{}{#1}{#2}}


\newenvironment{sarray}{
             \textfont0=\scriptfont0
             \scriptfont0=\scriptscriptfont0
             \textfont1=\scriptfont1
             \scriptfont1=\scriptscriptfont1
             \textfont2=\scriptfont2
             \scriptfont2=\scriptscriptfont2
             \textfont3=\scriptfont3
             \scriptfont3=\scriptscriptfont3
           \renewcommand{\arraystretch}{0.7}
           \begin{array}{l}}{\end{array}}

\newenvironment{scarray}{
             \textfont0=\scriptfont0
             \scriptfont0=\scriptscriptfont0
             \textfont1=\scriptfont1
             \scriptfont1=\scriptscriptfont1
             \textfont2=\scriptfont2
             \scriptfont2=\scriptscriptfont2
             \textfont3=\scriptfont3
             \scriptfont3=\scriptscriptfont3
           \renewcommand{\arraystretch}{0.7}
           \begin{array}{c}}{\end{array}}

\clearpage

\section{Introduction}  \label{sec1}

Let $\{c_n\}_{n=1}^\infty$ be a sequence of coefficients,
and let
\begin{equation}
    C(x) \;=\;  \sum_{n=1}^\infty c_n \, {x^n \over n!}
  \label{eq1.1}
\end{equation}
be the corresponding exponential generating function,
considered as a formal power series in the indeterminate $x$.
Now let $q$ be another (commuting) indeterminate\footnote{
   Let us stress that this $q$ has nothing to do with the $q$
   of $q$-series;
   rather, it is connected with the $q$ of the
   $q$-state Potts model \cite{Potts_52,Wu_82,Wu_84}
   and of the multivariate Tutte polynomial \cite{Sokal_bcc2005}.
},
and define
\begin{equation}
   A_q(x)  \;=\;  \exp[q C(x)]
           \;=\;  \sum_{n=0}^\infty a_n(q) \, {x^n \over n!}
   \;,
  \label{eq1.2}
\end{equation}
so that $A_q(x) = A(x)^q$ [we write $A(x) = A_1(x)$].
It is easy to see that:
\begin{itemize}
   \item[(a)]  Each coefficient $a_n(q)$ is a polynomial in $q$,
      of degree at most $n$; indeed we have the explicit formula\footnote{
   This formula holds also for $n=0$
   because the integer 0 has exactly one
   ordered partition $(n_1,\ldots,n_\ell)$ into
   zero or more strictly positive parts:
   namely, we have $\ell=0$ parts,
   and the multinomial coefficient and the empty product both equal 1,
   correctly giving $a_0(q) = 1$.
}
\begin{equation}
   a_n(q)  \;=\;
   \sum_{\ell=0}^n  {q^\ell \over \ell!}
   \sum_{\begin{scarray}
           n_1,\ldots,n_\ell \ge 1 \\
           \sum n_j = n
         \end{scarray}
        }
        {n \choose n_1,\ldots,n_\ell}
   \prod_{j=1}^\ell c_{n_j}
   \;.
   \label{eq.partition.an.cn}
\end{equation}
   \item[(b)]  $a_0(q)=1$.
   \item[(c)]  For $n \neq 0$, the polynomial $a_n(q)$ has zero constant term,
      so that we can write $a_n(q) = q \ahat_n(q)$
      where $\ahat_n(q)$ is a polynomial of degree at most $n-1$;
      moreover we have $\ahat_n(0) = c_n$.
\end{itemize}
For some purposes it is convenient to introduce the
modified generating function
\begin{equation}
   \widehat{A}_q(x)  \;=\;
  {A_q(x) \,-\, 1  \over q}  \;=\;
  {\exp[q C(x)] \,-\, 1  \over q}
           \;=\;  \sum_{n=1}^\infty \ahat_n(q) \, {x^n \over n!}
   \;.
  \label{eq1.2.hat}
\end{equation}

The coefficients $\{a_n(q)\}_{n=0}^\infty$ furthermore
satisfy the following identities \cite{Gould_74}:\footnote{
   See Gould \cite{Gould_74} for fascinating comments on the
   history of identities equivalent to \reff{eq.id2}
   [or to \reff{eq.recursion.anq1.anq2} below],
   which apparently go back to Euler (1748), Rothe (1793),
   Hindenburg (1796), von Ettingshausen (1826) and Hanstead (1881),
   among others.
}
\begin{eqnarray}
   \sum_{k=0}^n {n \choose k} \, a_k(q_1) \, a_{n-k}(q_2)
   & = &
   a_n(q_1 + q_2)
      \label{eq.id1} \\[2mm]
   \sum_{k=1}^n {n \choose k} \, k \, \ahat_k(q_1) \, a_{n-k}(q_2)
   & = &
   n \, \ahat_n(q_1 + q_2)
   \qquad\hbox{for } n \ge 1 \qquad
      \label{eq.id2}
\end{eqnarray}
Indeed, \reff{eq.id1} follows from writing the trivial identity
\begin{equation}
   A_{q_1}(x) A_{q_2}(x) \;=\; A_{q_1+q_2}(x)
\end{equation}
and extracting the coefficient of $x^n$,
while \reff{eq.id2} follows analogously from
\begin{equation}
   A_{q_2}(x) \left[ x {d \over dx} \, {A_{q_1}(x) \,-\, 1  \over   q_1}
              \right]
   \;\:=\;\:
   x {d \over dx} \, {A_{q_1+q_2}(x) \,-\, 1  \over   q_1 + q_2}
\end{equation}
\hspace*{-2mm}
\cite[eq.~(4.3)]{Gould_74},
which is an easy consequence of \reff{eq1.2}.
Specializing \reff{eq.id2} to $q_1 = 0$,
or alternatively applying
$x {\displaystyle d \over \displaystyle dx_{\vphantom{1}}}$
directly to \reff{eq1.2},
we obtain the well-known recursion relation \cite[Theorem 3.10.1]{Wilf_94}
\begin{equation}
   a_n(q)  \;=\;  q \sum_{k=1}^n  {n-1 \choose k-1} \, c_k \, a_{n-k}(q)
   \qquad\hbox{for } n \ge 1
   \;,
   \label{eq.recursion.an.cn}
\end{equation}
which allows the $\{a_n(q)\}$ to be calculated given the $\{c_n\}$,
or vice versa.
More generally, we can form the linear combination
$n \times \hbox{\reff{eq.id1}} - (q_1+q_2) \times \hbox{\reff{eq.id2}}$
to obtain the recursion relation\footnote{
   See \cite[eq.~(1.3)]{Gould_74}
   or \cite[Exercise~1.39, pp.~50 and 62]{Stanley_86}
   for the case $q_1 = 1$,
   which is of course trivially equivalent to the general case of $q_1 \neq 0$.
}
\begin{subeqnarray}
   a_n(q_2)
   & = &
   \sum_{k=1}^n
      \left[ {n-1 \choose k-1} (q_1+q_2) \,-\, {n \choose k} q_1 \right]
      \ahat_k(q_1) \, a_{n-k}(q_2)
   \\[2mm]
   & = &
   \sum_{k=1}^n
      \left[ {n-1 \choose k-1} q_2 \,-\, {n-1 \choose k} q_1 \right]
      \ahat_k(q_1) \, a_{n-k}(q_2)
   \label{eq.recursion.anq1.anq2}
\end{subeqnarray}
for $n \ge 1$,
which allows the $\{a_n(q_2)\}$ to be calculated given the $\{\ahat_n(q_1)\}$
for any fixed $q_1$,
and which reduces to \reff{eq.recursion.an.cn} when $q_1 = 0$.
We do not know whether there exist useful analogues of
\reff{eq.id1}/\reff{eq.id2} with higher powers of $k$
on the left-hand side.

Combinatorially, all these identities are versions of the
{\em exponential formula}\/
\cite{Uhlenbeck_62,Foata_70,Foata_74,Wilf_94,Bergeron_98,Stanley_99},
which relates the weights $c_n$ of ``connected'' objects
to the weights $a_n(q)$ of ``all'' objects,
when the weight of an object is defined to be the product of the weights
of its ``connected components'' times a factor $q$ for each
connected component.\footnote{
   In the {\em theory of combinatorial species}\/ \cite{Bergeron_98},
   the exponential formula is merely a special case of the
   more general concept of
   {\em composition}\/ (or {\em substitution}\/) {\em of species}\/ $F \circ G$
   \cite[section~1.4]{Bergeron_98},
   which is defined so that the exponential generating function of $F \circ G$
   is the composition of the exponential generating functions of $F$ and $G$.
   When $F$ is taken to be the species $E$ of sets
   ({\em ensembles}\/ in French)
   and $H$ is some species for which there exists a reasonable notion
   of connectedness,
   then, letting $H^c$ be the subspecies of connected $H$-structures,
   we have the isomorphism of species $H = E \circ H^c$;
   since the exponential generating function of $E$ is $\exp$,
   this immediately implies the exponential formula $H(x) = \exp(H^c(x))$.
   This result also generalizes to weighted and multi-sort species.
   See \cite{Labelle_96} for results closely related to those of the
   present paper, formulated in the language of species.
}
This idea can be formalized in the following combinatorial model:
Let $\Pi_n$ be the set of partitions of $[n] \equiv \{1,\ldots,n\}$
into (zero or more) nonempty blocks;
give a partition $\pi = \{\pi_1,\ldots,\pi_\ell\}$
a weight $w(\pi) = q^\ell \prod\limits_{i=1}^\ell c_{|\pi_i|}$,
where $\{c_k\}_{k=1}^\infty$ are arbitrary coefficients\footnote{
   The coefficients $c_k$ may in turn be generating functions
   for some additional structure on the set $[k]$
   --- for instance, a graph, a digraph, a hypergraph, etc.\  ---
   under the condition that this structure is ``connected''.
};
and let $a_n(q) = \sum\limits_{\pi \in \Pi_n} w(\pi)$.
Then \reff{eq.partition.an.cn} is immediate:
the multinomial coefficient counts the {\em ordered}\/ partitions of $[n]$
into blocks of sizes $n_1,\ldots,n_\ell$,
and the $1/\ell!$ converts this to {\em unordered}\/ partitions.
As for \reff{eq.id1}/\reff{eq.id2},
we first observe that $a_n(q_1+q_2)$ can be computed by allowing each block
to be colored either ``1'' (with weight $q_1$) or ``2'' (with weight $q_2$).
Identity \reff{eq.id1} follows by defining $S \subseteq [n]$
to be the set of elements colored ``1'', and setting $k = |S|$.
Analogously, $\ahat_n(q_1+q_2)$ can be computed by coloring and weighting
blocks as before, except that the block containing one specified
element $i \in [n]$ (say, $i=1$) must be colored ``1''
and it receives weight 1 instead of weight $q_1$.
Defining again $S \subseteq [n]$ to be the set of elements colored ``1'',
and setting $k = |S|$,  we obtain \reff{eq.id2} in the form\footnote{
   See also \cite[Theorem~2.9]{Reiner_77};
   and see \cite[p.~188, Exercise~3.1.20(a)]{Bergeron_98}
   for the special case $q_1 = 1$.
}
\begin{equation}
   \ahat_n(q_1 + q_2)
   \;=\;
   \sum_{k=1}^n {n-1 \choose k-1} \, \ahat_k(q_1) \, a_{n-k}(q_2)
   \qquad\hbox{for } n \ge 1
   \;.
 \label{eq.id2.bis}
\end{equation}


Sequences of polynomials $\{a_n(q)\}$ satisfying \reff{eq1.2}
[or equivalently satisfying \reff{eq.id1} and not identically zero]
have been termed {\em sequences of binomial type}\/
by Rota and collaborators
\cite{Mullin_70,Rota_73,Garsia_73,Fillmore_73,Roman_78,Roman_84}
and studied by means of the umbral calculus
\cite{Roman_84,Rota_94,Gessel_03}.\footnote{
   It is not hard to see by induction that if $\{a_n(q)\}$
   satisfies \reff{eq.id1}
   [or even the special case of \reff{eq.id1} with $q_1 = q_2$]
   and is not identically zero, then $\deg a_n \le n$ with $a_0 = 1$
   and $a_n(0) = 0$ for $n \ge 1$;
   moreover, with a little extra work it can be shown
   that $A_q(x)$ is necessarily of the form \reff{eq1.2}:
   see \cite[Theorem~4.1]{Garsia_73} \cite{Knuth_92,Zeng_96}
   \cite[Exercise~5.37, pp.~87--88 and 131--132]{Stanley_99}.

   Let us also remark that Rota {\em et al.}\/ \cite{Mullin_70,Rota_73}
   and many subsequent authors \cite{Garsia_73,Fillmore_73,Roman_78,Roman_84},
   in defining ``sequence of binomial type'',
   impose the additional condition that $\deg a_n = n$ exactly
   (i.e.\ $c_1 \neq 0$);
   but since this condition is irrelevant for our purposes,
   we prefer not to impose it.
}
The corresponding sequences $\{a_n(q)/n!\}$
have been termed sequences of binomial type by Labelle \cite{Labelle_80}
and {\em convolution families}\/
by Knuth \cite{Knuth_92} and Zeng \cite{Zeng_96};
these authors used elementary formal-power-series methods
closely resembling those used here.
A purely combinatorial approach to sequences of binomial type,
employing the theory of species,
has been developed by Labelle \cite{Labelle_81}
(see also \cite[section~3.1]{Bergeron_98}).

Further identities, generalizing \reff{eq.id1}/\reff{eq.id2},
can be obtained by a powerful transformation
suggested by Knuth \cite{Knuth_92}.
Start with any formal power series $A(x)$ with constant term 1;
set $A_q(x) = A(x)^q$ and define polynomials $\{a_n(q)\}$ by \reff{eq1.2}.
Now let $t$ be a parameter, and define a formal power series $A(x;t)$
by the implicit equation
\begin{equation}
   A(x;t)   \;=\;  A\bigl( x A(x;t)^t \bigr)
   \;.
 \label{eq.Axt.knuth}
\end{equation}
Since $A(x;0) = A(x)$,
we can view the family $\{A(x;t)\}_{t \in \R}$
as a special class of ``perturbations'' of $A(x)$.
Now define polynomials $\{a_n(q;t)\}$ by the obvious formula
\begin{equation}
   A(x;t)^q   \;=\;  \sum_{n=0}^\infty a_n(q;t) \, {x^n \over n!}
   \;.
  \label{eq1.2.t}
\end{equation}
For each $t$, these polynomials form a family of the type \reff{eq1.2}
and thus satisfy the identities \reff{eq.id1}/\reff{eq.id2}.
On the other hand,
a straightforward calculation using the Lagrange inversion formula
\cite{Stanley_99} yields the remarkable relation\footnote{
   See \cite[Exercises~5.37(e) and 5.58,
             pp.~87--88, 99, 131--133 and 148]{Stanley_99},
   where the argument is attributed to Eric Rains and Ira Gessel.
   Please note also that if we write
   $f(x) = x/A(x)^t$ and $g(x) = x A(x;t)^t$,
   then $f$ and $g$ are compositional inverses
   (as observed in \cite{Labelle_80,Knuth_92} for $t=1$).
}
\begin{equation}
   a_n(q;t)  \;=\;  {q \over q+nt} \: a_n(q+nt)   \;=\;  q \, \ahat_n(q+nt)
   \;,
 \label{eq.anqt}
\end{equation}
or even more simply
\begin{equation}
   \ahat_n(q;t)   \;=\;   \ahat_n(q+nt)
   \;.
\end{equation}
[In particular, $a_n(q;t)$ and $\ahat_n(q;t)$ are {\em polynomials}\/
 jointly in $q$ and $t$, so we can treat $t$ as an indeterminate if we wish.]
It follows that, for any family $\{a_n(q)\}$ of type \reff{eq1.2},
we have the following identities generalizing
\reff{eq.id1}/\reff{eq.id2}:\footnote{
   Already Rota {\em et al.}\/ \cite[p.~711, Proposition~4]{Rota_73}
   and Reiner \cite[Theorem 2.10]{Reiner_77}
   noticed that, for any $t$, the family $\{a_n(q;t)\}$
   defined by \reff{eq.anqt} is of type \reff{eq1.2}.
   But they do not seem to have noticed the relation \reff{eq.Axt.knuth},
   which to our knowledge goes back to Knuth \cite{Knuth_92}.
}
\begin{eqnarray}
   \sum_{k=0}^n {n \choose k} \, q_1 \ahat_k(q_1+kt)
                              \, q_2 \ahat_{n-k}(q_2+(n-k)t)
   & = &
   (q_1+q_2) \ahat_n(q_1 + q_2+nt)
      \label{eq.id1.t} \\[2mm]
   \sum_{k=1}^n {n \choose k} \, k \, \ahat_k(q_1+kt)
                                   \, q_2 \ahat_{n-k}(q_2+(n-k)t)
   & = &
   n \, \ahat_n(q_1 + q_2+nt)
   \qquad\hbox{for } n \ge 1
      \nonumber \\[-4mm] \label{eq.id2.t}
\end{eqnarray}
where we understand that $q \ahat_0(q) = 1$.
Multiplying \reff{eq.id2.t} by $t$ and summing it with \reff{eq.id1.t},
we obtain the alternate form
\begin{equation}
   \sum_{k=0}^n {n \choose k} \, a_k(q_1+kt)
                              \, q_2 \ahat_{n-k}(q_2+(n-k)t)
   \;=\;
   a_n(q_1 + q_2+nt)
   \;.
      \label{eq.id3.t}
\end{equation}
If in \reff{eq.id3.t} we interchange $k \leftrightarrow n-k$
and $q_1 \leftrightarrow q_2$ and then replace $q_2$ by $q_2 -nt$,
we obtain the slightly simpler form
\begin{equation}
   \sum_{k=0}^n {n \choose k} \, q_1 \ahat_k(q_1+kt)
                              \, a_{n-k}(q_2-kt)
   \;=\;
   a_n(q_1 + q_2)
   \;.
      \label{eq.id4.t}
\end{equation}
A combinatorial proof and interpretation of \reff{eq.id1.t}/\reff{eq.id3.t}
in terms of the theory of species
has been given by Labelle \cite{Labelle_81}.\footnote{
   See also \cite[p.~189, Exercise~3.1.22]{Bergeron_98} for a brief summary.
}
We call \reff{eq.id1.t}--\reff{eq.id4.t}
{\em Abel-type extensions}\/ of \reff{eq.id1}/\reff{eq.id2}
because in the simplest case $A(x) = e^x$, $a_n(q) = q^n$
they reduce to Abel's celebrated extensions of the binomial theorem
(also discovered simultaneously by Cauchy)
\cite{Riordan_68,Strehl_92,Johnson_07}.

These Abel-type formulae can alternatively be derived
as immediate consequences of the Pfaff--Cauchy derivative identities
(generalizations of Leibniz's rule for the $n$th derivative of the
product of two functions) recently presented by Johnson \cite{Johnson_07}.
Indeed, setting $v = A(x)^{q_1}$, $w = A(x)^{q_2}$, $\phi = A(x)^t$
in \cite[eq.~(1.1)]{Johnson_07} yields \reff{eq.id4.t},
while making the same substitution in
\cite[eq.~(1.4)]{Johnson_07} yields \reff{eq.id1.t}.

Some classical examples of families $\{a_n(q)\}$
and the corresponding identities
\reff{eq.id1}/ \reff{eq.id2}
and \reff{eq.id1.t}--\reff{eq.id4.t}
are listed in Table~\ref{table1}.

\begin{table}[t]
\addtolength{\parskip}{-3mm}
\hspace*{-1.2cm}
\begin{tabular}{|c|c|c|c|}
\hline
  $A(x)$   &   $a_n(q)$   &   \reff{eq.id1}/\reff{eq.id2}   &
                              \reff{eq.id1.t}--\reff{eq.id4.t}   \\
\hline
   $e^x$   &   $q^n$
           &   $\begin{array}{c}
                   \hbox{\small binomial} \\[-1mm] \hbox{\small theorem}
                \end{array}$
           &   $\begin{array}{c}
                   \hbox{\small Abel--Cauchy} \\[-1mm]
                   \hbox{\small binomial} \\[-1mm] \hbox{\small identities}
                \end{array}$
           \\[3mm]
   $\begin{array}{c}
       1+x          \\[3mm]
       (1-x)^{-1}   \\[3mm]
       (1+\alpha x)^\beta
    \end{array}$
   &
   $\left.
    \begin{array}{c}
       q^{\underline{n}}   \\[3mm]
       q^{\overline{n}}    \\[3mm]
       \alpha^n (\beta q)^{\underline{n}}
    \end{array}
    \!\right\}$
           &   $\begin{array}{c}
                   \hbox{\small Chu--Vandermonde} \\[-1mm]
                   \hbox{\small convolution}
                \end{array}$
           &   $\begin{array}{c}
                   \hbox{\small Rothe's generalized} \\[-1mm]
                   \hbox{\small Chu--Vandermonde} \\[-1mm]
                   \hbox{\small convolution}
                \end{array}$
          \\
   &  &  &  \\[-1mm]
   $\exp(e^x - 1)$
   &
   $\begin{array}{c}
      B_n(q) \,=\,
      \sum\limits_{k=0}^n \stirlingsubset{n}{k} \, q^k \\[1mm]
      \hbox{\small (Bell polynomial \cite{Comtet_74})}
   \end{array}$
          &  
          &  
          \\[9mm]
   $\exp[-x/(1-x)]$
   &
   $\begin{array}{c}
      L_n^{(-1)}(q)  \,=\,
      \sum\limits_{k=0}^n {n! \over k!} \, {n-1 \choose n-k} \, (-q)^k \\[1mm]
      \hbox{\small (Laguerre polynomial \cite{Andrews_99})}
   \end{array}$
          &  
          &  
          \\
\hline
\end{tabular}
\caption{
   Some classical identities that are special cases of
   \reff{eq.id1}/\reff{eq.id2} and
   \reff{eq.id1.t}--\reff{eq.id4.t}.
   See \cite{Gould_56,Gould_57,Gould_66,Gould_74,Riordan_68,Strehl_92,%
Johnson_07}
   for more information concerning Abel--Cauchy and
   Chu--Vandermonde--Rothe identities.
   The identity \reff{eq.id1} is well known for the Bell
   \cite[p.~136, eq.~3n]{Comtet_74}
   and Laguerre
   \cite[p.~192, eq.~41]{Bateman} \cite[p.~96, Problem~20]{Lebedev_65}
   polynomials.
   We do not know any references for \reff{eq.id2} or 
   \reff{eq.id1.t}--\reff{eq.id4.t} in these cases.
}
\label{table1}
\addtolength{\parskip}{3mm}
\end{table}

\bigskip

In this note we begin by proving some easy generalizations
of the foregoing formulae to the case in which the integer index $n$
is replaced by a multi-index $\n$ over a finite ground set $V$.
We then show some applications of the latter formulae
to the multivariate Tutte polynomial \cite{Sokal_bcc2005}
of a graph $G=(V,E)$.
Finally, we derive some further identities for the
multivariate Tutte polynomial,
which generalize an identity for counting connected graphs
found by Riordan, Nijenhuis, Wilf and Kreweras
\cite{Riordan_unpub,Nijenhuis_79,Kreweras_80}
and in more general form
by Leroux \cite{Leroux_88} and Gessel \cite{Gessel_95}
and an identity for the inversion enumerator of trees
found by Mallows, Riordan and Kreweras \cite{Mallows_68,Kreweras_80}.
In the appendix we prove an apparently new generalization
of M\"obius inversion on the partition lattice.

We stress that many (though not all) of the formulae
derived in this paper are already known.
But we think that our approach provides a unified perspective
that may be of interest.

\section{Variants of the exponential formula}  \label{sec2}

Let $V$ be a finite ground set;
we shall use multi-indices $\n = (n_i)_{i \in V} \in \N^V$
and commuting indeterminates $\x = (x_i)_{i \in V}$.
We write $|\n| = \sum_{i \in V} n_i$
and $\x^\n = \prod_{i \in V} x_i^{n_i}$.
We use the notation $\n! = \prod_{i \in V} n_i!$
and analogously for binomial and multinomial coefficients.
Finally, we denote by $\zero$ (resp.\ $\one$)
the multi-index with all entries 0 (resp.\ 1);
and for any subset $W \subseteq V$, we denote by $\one_W$
the multi-index taking the value 1 on $W$ and 0 on $V \setminus W$.

Let $R$ be a commutative ring containing the rationals as a subring;
we shall use polynomials and formal power series
whose coefficients lie in $R$.
(In applications, $R$ will usually be either the rationals,
 reals or complex numbers, or a ring of polynomials or formal power series
 over one of these fields.)
One possible approach is to treat the coefficients
$\cc = (c_\n)_{\n \in \N^V \setminus \{\zero\}}$
as indeterminates, and to take $R$ to be the polynomial ring $\Q[\cc]$;
this approach has the advantage of exhibiting the polynomial dependence
on the coefficients $c_\n$.
But we need not commit ourselves to any specific choice of the ring $R$;
all of our identities will be valid in complete generality.

So let $\{c_\n\}_{\n \in \N^V}$ be a sequence of coefficients in the ring $R$,
with $c_\zero = 0$;
and let
\begin{equation}
    C(\x) \;=\;  \sum_{\n} c_\n \, {\x^\n \over \n!}
\end{equation}
be the corresponding exponential generating function,
considered as a formal power series in the indeterminates $\x$.
Now let $q$ be another (commuting) indeterminate, and define
\begin{equation}
   A_q(\x)  \;=\;  \exp[q C(\x)]
           \;=\;  \sum_{\n} a_\n(q) \, {\x^\n \over \n!}
   \;.
 \label{def.Aq.multidim}
\end{equation}
Using the Taylor series for $\exp$, we deduce immediately that
\begin{itemize}
   \item[(a)]  Each coefficient $a_\n(q)$ is a polynomial in $q$,
      of degree at most $|\n|$;
      indeed we have the explicit formula\footnote{
   This formula holds also for $\n = \zero$
   because the multi-index $\zero$ has exactly one
   ordered partition $(\n_1,\ldots,\n_\ell)$ into
   zero or more nonvanishing parts:
   namely, we have $\ell=0$ parts,
   and the multinomial coefficient and the empty product both equal 1,
   correctly giving $a_\zero(q) = 1$.
   This same observation applies to all subsequent formulae
   of a similar type, such as \reff{eq.partition2.rq.q} ff.
}
\begin{equation}
   a_\n(q)  \;=\;
   \sum_{\ell=0}^{|\n|}  {q^\ell \over \ell!}
   \sum_{\substack{
           \n_1,\ldots,\n_\ell \neq \zero \\[0.7mm]
           \sum \n_j = \n
         }
        }
        {\n \choose \n_1,\ldots,\n_\ell}
   \prod_{j=1}^\ell c_{\n_j}
   \;.
   \label{eq.partition2.an.cn}
\end{equation}
   \item[(b)]  $a_\zero(q) = 1$.
   \item[(c)]  For $\n \neq \zero$, the polynomial $a_\n(q)$
      has zero constant term,
      so that we can write $a_\n(q) = q \ahat_\n(q)$
      where $\ahat_\n(q)$ is a polynomial of degree at most $|\n|-1$;
      moreover we have $\ahat_\n(0) = c_\n$.
\end{itemize}
Conversely, it is not hard to see that
if $(c_\n)_{\n \in \N^V \setminus \{\zero\}}$ is any sequence
and we define $(a_\n(q))_{\n \in \N^V}$ by \reff{eq.partition2.an.cn},
then the latter family satisfies \reff{def.Aq.multidim}.

Using $q C(\x) = \log A_q(\x)$
together with the Taylor series for $\log(1+z)$,
we obtain the inverse formula
\begin{equation}
   q c_\n  \;=\;
   \sum_{\ell=1}^{|\n|}  {(-1)^{\ell-1} \over \ell}
   \sum_{\substack{
           \n_1,\ldots,\n_\ell \neq \zero \\[0.7mm]
           \sum \n_j = \n
         }
        }
        {\n \choose \n_1,\ldots,\n_\ell}
   \prod_{j=1}^\ell a_{\n_j}(q)
 \label{eq.partition2.inverse.1}
\end{equation}
or equivalently
\begin{equation}
   c_\n  \;=\;
   \sum_{\ell=1}^{|\n|}  {(-q)^{\ell-1} \over \ell}
   \sum_{\substack{
           \n_1,\ldots,\n_\ell \neq \zero \\[0.7mm]
           \sum \n_j = \n
         }
        }
        {\n \choose \n_1,\ldots,\n_\ell}
   \prod_{j=1}^\ell \ahat_{\n_j}(q)
   \;.
 \label{eq.partition2.inverse.2}
\end{equation}
It is a nontrivial fact that the right-hand side of
 \reff{eq.partition2.inverse.2} is independent of $q$.

Now let $r$ be another indeterminate,
and let us use the identity $A_{rq}(\x) = A_q(\x)^r$
together with the Taylor series
\begin{equation}
   (1+z)^r  \;\equiv\; \exp[r \log(1+z)]  \;=\;
   \sum_{\ell=0}^\infty  {r^{\underline{\ell}}  \over  \ell!} \, z^\ell
 \label{eq.Taylor.power_r}
\end{equation}
where $r^{\underline{\ell}} = r (r-1) \cdots (r-\ell+1)$
denotes the falling factorial.
(If we like, we can alternatively use the notation
 ${r \choose \ell} = r^{\underline{\ell}}/\ell!$,
 which is a well-defined polynomial in the indeterminate $r$.)
We obtain
\begin{equation}
   a_\n(rq)  \;=\;
   \sum_{\ell=0}^{|\n|}  {r^{\underline{\ell}}  \over  \ell!}
   \sum_{\substack{
           \n_1,\ldots,\n_\ell \neq \zero \\[0.7mm]
           \sum \n_j = \n
         }
        }
        {\n \choose \n_1,\ldots,\n_\ell}
   \prod_{j=1}^\ell a_{\n_j}(q)
   \;.
 \label{eq.partition2.rq.q}
\end{equation}
In view of the identity
$(-r)^{\underline{\ell}} = (-1)^\ell r^{\overline{\ell}}$
where $r^{\overline{\ell}} = r (r+1) \cdots (r+\ell-1) =
      (r+\ell-1)^{\underline{\ell}}$
denotes the rising factorial,
we also have
\begin{equation}
   a_\n(-rq)  \;=\;
   \sum_{\ell=0}^{|\n|} {(-1)^\ell r^{\overline{\ell}}  \over  \ell!}
   \sum_{\substack{
           \n_1,\ldots,\n_\ell \neq \zero \\[0.7mm]
           \sum \n_j = \n
         }
        }
        {\n \choose \n_1,\ldots,\n_\ell}
   \prod_{j=1}^\ell a_{\n_j}(q)
   \;.
 \label{eq.partition2.minusrq.q}
\end{equation}

Since both sides of the identities
\reff{eq.partition2.rq.q}/\reff{eq.partition2.minusrq.q}
are polynomials in $r$,
we can specialize these identities to $r$ integer or rational,
or more generally to $r$ being any element of the ring $R$.
If $r$ is a positive integer, the sum in \reff{eq.partition2.rq.q}
can obviously be restricted to $\ell \le r$.
Slightly (but not much) less trivially,
if $r$ is a positive integer, then by using
\begin{equation}
   {r^{\overline{\ell}}  \over  \ell!}
   \;=\;
   {(r+\ell-1)!  \over \ell! \, (r-1)!}
   \;=\;
   {(\ell+1)^{\overline{r-1}} \over (r-1)!}
\end{equation}
we can rewrite \reff{eq.partition2.minusrq.q} as
\begin{equation}
   a_\n(-rq)  \;=\;
   \sum_{\ell=0}^{|\n|} (-1)^\ell {(\ell+1)^{\overline{r-1}} \over (r-1)!}
   \sum_{\substack{
           \n_1,\ldots,\n_\ell \neq \zero \\[0.7mm]
           \sum \n_j = \n
         }
        }
        {\n \choose \n_1,\ldots,\n_\ell}
   \prod_{j=1}^\ell a_{\n_j}(q)
   \;,
 \label{eq.partition2.minusrq.q.bis}
\end{equation}
in which the coefficient $(\ell+1)^{\overline{r-1}} / (r-1)!$
is a polynomial of degree $r-1$ in $\ell$.
In particular, for $r=1$ we have
\begin{equation}
   a_\n(-q)  \;=\;
   \sum_{\ell=0}^{|\n|}  (-1)^{\ell}
   \sum_{\substack{
           \n_1,\ldots,\n_\ell \neq \zero \\[0.7mm]
           \sum \n_j = \n
         }
        }
        {\n \choose \n_1,\ldots,\n_\ell}
   \prod_{j=1}^\ell a_{\n_j}(q)
   \;.
   \label{eq.partition2bis.minusq}
\end{equation}


By formally setting $r=q_2/q_1$ (this is justifiable if we work in a ring
of formal Laurent series in the indeterminate $q_1$),
we obtain the identity
\begin{equation}
   a_\n(q_2)  \;=\;
   \sum_{\ell=0}^{|\n|}  {1 \over \ell!}
   \, \Biggl( \, \prod_{j=0}^{\ell-1} (q_2 - jq_1) \, \Biggr)
   \sum_{\substack{
           \n_1,\ldots,\n_\ell \neq \zero \\[0.7mm]
           \sum \n_j = \n
         }
        }
        {\n \choose \n_1,\ldots,\n_\ell}
   \prod_{j=1}^\ell \ahat_{\n_j}(q_1)
   \;,
   \label{eq.partition2bis.q1.q2}
\end{equation}
which is nice because both sides are manifestly polynomials in the
indeterminates $q_1$ and $q_2$.
Alternatively, for $\n \neq \zero$ we can pull out the $j=0$ factor
and write
\begin{equation}
   \ahat_\n(q_2)  \;=\;
   \sum_{\ell=1}^{|\n|}  {1 \over \ell!}
   \, \Biggl( \, \prod_{j=1}^{\ell-1} (q_2 - jq_1) \, \Biggr)
   \sum_{\substack{
           \n_1,\ldots,\n_\ell \neq \zero \\[0.7mm]
           \sum \n_j = \n
         }
        }
        {\n \choose \n_1,\ldots,\n_\ell}
   \prod_{j=1}^\ell \ahat_{\n_j}(q_1)
   \;.
   \label{eq.partition2bis.q1hat.q2}
\end{equation}
Conversely, it is not hard to see that
if $q_1$ is any element of the ring $R$
and $(\ahat_\n(q_1))_{\n \in \N^V}$ is any sequence,
and we define $(a_\n(q))_{\n \in \N^V}$ by \reff{eq.partition2bis.q1.q2},
then the latter family satisfies \reff{def.Aq.multidim}.

The following special cases of
\reff{eq.partition2bis.q1.q2}/\reff{eq.partition2bis.q1hat.q2}
are of particular interest:
\begin{itemize}
   \item[(a)]  For $q_1=0$, the identity \reff{eq.partition2bis.q1.q2}
       reduces to \reff{eq.partition2.an.cn}.
   \item[(b)]  For $q_2=0$, we obtain the inverse formula
       \reff{eq.partition2.inverse.1}/\reff{eq.partition2.inverse.2}.
   \item[(c)] For $q_2 = -q_1$, we obtain the closely related formula
       \reff{eq.partition2bis.minusq}.
\end{itemize}

It is easy to show that the coefficients $\{a_\n(q)\}_{\n \in \N^V}$
satisfy the identities
\begin{eqnarray}
   \sum_{\k} {\n \choose \k} \, a_\k(q_1) \, a_{\n-\k}(q_2)
   & = &
   a_\n(q_1 + q_2)
      \label{eq.id1.multidim}  \\[2mm]
   \sum_{\k \neq \zero}
      {\n \choose \k} \, k_i \, \ahat_\k(q_1) \, a_{\n-\k}(q_2)
   & = &
   n_i \, \ahat_\n(q_1 + q_2)
   \qquad\hbox{for } \n \neq \zero \qquad
      \label{eq.id2.multidim}
\end{eqnarray}
Indeed, \reff{eq.id1.multidim} follows by expanding out the identity
$A_{q_1}(\x) A_{q_2}(\x) = A_{q_1+q_2}(\x)$
and extracting the coefficient of $\x^\n$,
while \reff{eq.id2.multidim} follows analogously from
\begin{equation}
   A_{q_2}(\x) \left[ x_i {\partial  \over \partial x_i} \,
                      {A_{q_1}(\x) \,-\, 1  \over   q_1}
               \right]
   \;\:=\;\:
   x_i {\partial  \over \partial x_i} \,
       {A_{q_1+q_2}(\x) \,-\, 1  \over   q_1 + q_2}
   \;,
\end{equation}
which is an easy consequence of \reff{def.Aq.multidim}.
In particular, specializing \reff{eq.id2.multidim} to $q_1 = 0$,
we obtain the recursion relation
\begin{equation}
   a_\n(q)  \;=\; q \sum_{\k \ge \bdelta_i}
       {\n-\bdelta_i \choose \k-\bdelta_i}  \, c_\k \, a_{\n-\k}(q)
   \qquad\hbox{whenever } n_i \ge 1 \;, \qquad
   \label{eq.recursion.an.cn.multidim}
\end{equation}
which allows the $\{a_\n(q)\}$ to be calculated given the $\{c_\n\}$,
or vice versa.
[Here $\bdelta_i$ is the vector with entry 1 at element $i$
 and 0 elsewhere.]
More generally, we can form the linear combination
$n_i \times \hbox{\reff{eq.id1.multidim}} -
 (q_1+q_2) \times \hbox{\reff{eq.id2.multidim}}$
to obtain the recursion relation
\begin{equation}
   a_\n(q_2)  \;=\;
   \sum_{\k \neq \zero}
     \left[ {\n-\bdelta_i \choose \k-\bdelta_i} (q_1+q_2)
            \,-\, 
            {\n \choose \k} q_1
     \right]
      \ahat_\k(q_1) \, a_{\n-\k}(q_2)
   \qquad\hbox{whenever } n_i \ge 1 \;, \qquad
   \label{eq.recursion.anq1.anq2.multidim}
\end{equation}
which allows the $\{a_\n(q_2)\}$ to be calculated given the $\{\ahat_\n(q_1)\}$
for any fixed $q_1$,
and which reduces to \reff{eq.recursion.an.cn.multidim} when $q_1 = 0$.

We do not know whether there exist useful analogues of
\reff{eq.id1.multidim}/\reff{eq.id2.multidim} with higher powers of $\k$
on the left-hand side.

Finally, we can deduce the Abel-type extensions of
\reff{eq.id1.multidim}/\reff{eq.id2.multidim}:
\begin{eqnarray}
   \sum_{\k} {\n \choose \k} \, q_1 \ahat_\k(q_1+ \k\cdot\ttt)
                              \, q_2 \ahat_{\n-\k}(q_2+(\n-\k)\cdot\ttt)
   & = &
   (q_1+q_2) \ahat_\n(q_1 + q_2+\n\cdot\ttt)
      \nonumber \\[-4mm] \label{eq.id1.t.multidim} \\[2mm]
   \sum_{\k \neq \zero} {\n \choose \k} \, k_i \, \ahat_\k(q_1+\k\cdot\ttt)
                                   \, q_2 \ahat_{\n-\k}(q_2+(\n-\k)\cdot\ttt)
   & = &
   n_i \, \ahat_\n(q_1 + q_2+\n\cdot\ttt)
   \quad\hbox{for } \n \neq \zero  \hspace*{-1cm}
      \nonumber \\[-4mm] \label{eq.id2.t.multidim}
\end{eqnarray}
where $\ttt = (t_i)_{i \in V}$ are commuting indeterminates
and we understand $q \ahat_\zero(q) = 1$.
These formulae can be proved by repeated use of the Lagrange inversion formula,
once in each variable.
For details we refer to the paper of Zeng \cite[pp.~225--226]{Zeng_96}.
For further information on classical specializations of
\reff{eq.id1.t.multidim}/\reff{eq.id2.t.multidim},
see \cite{Stam_87,Strehl_92,Zeng_96,Pitman_01,Pitman_02}.

\medskip

{\bf Remark.}  We wonder whether
\reff{eq.id1.t.multidim}/\reff{eq.id2.t.multidim}
might have an extension to identities involving a
not-necessarily-diagonal {\em matrix}\/ $T = (t_{ij})_{i,j \in V}$
of indeterminates, based on considering the implicit equation
\begin{equation}
   e^{C({\bf x};T)} = A\bigl( e^{C({\bf x};T) T} {\bf x} \bigr)
   \;,
\end{equation}
which generalizes \reff{eq.Axt.knuth}.
If true, this could be a powerful extension.
See Strehl \cite{Strehl_92} for possibly related work
in which a matrix of indeterminates is employed.

\section{Application to the multivariate Tutte polynomial}  \label{sec3}

In this section we apply the formulae of Section~\ref{sec2}
to deduce some identities for multivariate Tutte polynomials.

\subsection{Definitions and basic properties}

Let $G = (V,E)$ be a finite undirected graph
with vertex set $V \neq \emptyset$ and edge set $E$;
loops and multiple edges are allowed unless explicitly stated otherwise.
Then the {\em multivariate Tutte polynomial}\/ \cite{Sokal_bcc2005}
of $G$ is, by definition, the polynomial
\begin{equation}
   Z_G(q, \bv)   \;=\;
   \sum_{A \subseteq E}  q^{k(A)}  \prod_{e \in A}  v_e
   \;,
   \label{def.ZG}
\end{equation}
where $q$ and $\bv = (v_e)_{e \in E}$ are commuting indeterminates,
and $k(A)$ denotes the number of connected components in the subgraph $(V,A)$.
It is also convenient to pull out one factor of $q$ by defining
\begin{equation}
   \Zhat_G(q, \bv)   \;=\;
   \sum_{A \subseteq E}  q^{k(A)-1}  \prod_{e \in A}  v_e
   \;;
   \label{def.ZhatG}
\end{equation}
this is still a polynomial since $k(A) \ge 1$ for all $A$.
Finally, let us make the convention that if $G=\emptyset$
(the graph with empty vertex set and empty edge set),
then $Z_\emptyset = 1$ and $\Zhat_\emptyset$ is undefined.

If we specialize to $v_e = v$ for all edges $e$,
we obtain a two-variable polynomial $Z_G(q,v)$
that is essentially equivalent to the classical Tutte polynomial
\cite[section~2.5]{Sokal_bcc2005}.

Note that at $q=1$ we have the trivial formula
\begin{equation}
   Z_G(1, \bv)  \;=\;  \Zhat_G(1, \bv)  \;=\;   \prod_{e \in E} (1+v_e)  \;.
 \label{ZG.q=1}
\end{equation}
More interestingly, at $q=0$ we have the specialization
\begin{equation}
   \Zhat_G(0, \bv)  \;=\;
   C_G(\bv) \;\equiv\; \!\! \sum\limits_{\begin{scarray}
                                         A \subseteq E \\
                                         k(A) = 1
                                     \end{scarray}}
                        \prod_{e \in A}  v_e
   \;,
\end{equation}
i.e.\ the generating polynomial of connected spanning subgraphs of $G$.

In statistical physics, $Z_G(q,\bv)$ is known as the partition function
of the $q$-state Potts model \cite{Potts_52,Wu_82,Wu_84}
in the Fortuin--Kasteleyn representation
\cite{Kasteleyn_69,Fortuin_72,Grimmett_book_to_appear}.
This arises by virtue of the following identity
\cite{Kasteleyn_69,Fortuin_72}:

\begin{theorem}[Fortuin--Kasteleyn representation of the Potts model
                   \protect\cite{Kasteleyn_69,Fortuin_72}]
   \label{thm.FK}
For integer $q \ge 1$, we have
\begin{equation}
   Z_G(q, \bv)  \;=\;
   \sum_{ \sigma \colon\, V \to [q]}
   \; \prod_{e=ij \in E}  \,
      \biggl[ 1 + v_e \delta\bigl(\sigma(i), \sigma(j)\bigr) \biggr]
 \label{eq.FK.identity}
\end{equation}
where $\delta$ denotes the Kronecker delta.
\end{theorem}

\proof
On the right-hand side of \reff{eq.FK.identity},
expand out the product over $e \in E$,
and let $A \subseteq E$ be the set of edges for which the term
$v_e \delta\bigl(\sigma(i), \sigma(j)\bigr)$ is taken.
Now perform the sum over maps $\sigma \colon\, V \to [q]$:
in each connected component of the subgraph $(V,A)$
the ``color'' $\sigma(i)$ must be constant,
and there are no other constraints.
Therefore, the right-hand side equals
\begin{equation}
   \sum_{ A \subseteq E }  q^{k(A)}  \prod_{e \in A}  v_e
   \;,
\end{equation}
as was to be proved.
\qed

Specializing to $v_e =-1$ for all edges $e$,
we obtain the Birkhoff--Whitney \cite{Birkhoff_12,Whitney_32a}
expansion for the chromatic polynomial of $G$:

\begin{corollary}
  \label{cor.chrompoly}
For integer $q \ge 1$, the number of proper $q$-colorings of $G$
is $P_G(q) \equiv Z_G(q,-1)$.
\end{corollary}

\subsection{Some generating functions}

For the remainder of this section,
let us assume that $G$ is a {\em loopless}\/ graph.
For notational simplicity it is convenient to assume also
that $G$ has no multiple edges;
the trivial changes to allow multiple edges can be left to the reader.

For any multi-index $\n \in \N^V$, let us define $G[\n]$
to be the graph obtained from $G$ by expanding each vertex $i$
to an independent set consisting of $n_i$ vertices.
That is, the vertices of $G[\n]$ are pairs $(i,\alpha)$
with $i \in V$ and $\alpha \in [n_i]$,
and the edges of $G[\n]$ are all pairs $\< (i,\alpha), (j,\beta) \>$
with $ij \in E$, $\alpha \in [n_i]$ and $\beta \in [n_j]$.
Given a set of weights $\bv = (v_e)_{e \in E}$
associated to the edges of $G$,
we assign the weight $v_{ij}$ to each edge
$\< (i,\alpha), (j,\beta) \>$ in $G[\n]$.

In a similar way, let us define $G'[\n]$
to be the graph obtained from $G$ by expanding each vertex $i$
to a clique consisting of $n_i$ vertices.
That is, the vertex set of $G'[\n]$ is the same as that of $G[\n]$,
and the edges of $G'[\n]$ consist of those of $G[\n]$
together with all pairs $\< (i,\alpha), (i,\beta) \>$
with $i \in V$, $\alpha,\beta \in [n_i]$ and $\alpha \neq \beta$.
Given a set of weights $\bv = (v_e)_{e \in E}$
associated to the edges of $G$
and another set of weights $\w = (w_i)_{i \in V}$
associated to the vertices of $G$,
we assign the weight $v_{ij}$ to each edge
$\< (i,\alpha), (j,\beta) \>$ ($i \neq j$) in $G'[\n]$
and the weight $w_i$ to each edge
$\< (i,\alpha), (i,\beta) \>$ in $G'[\n]$.

Note in particular that if $\n = \one_W$ for some subset $W \subseteq V$,
then $G[\n] = G'[\n] = G[W]$, the induced subgraph of $G$ on $W$.

The following ``master formula''
generates the multivariate Tutte polynomials
(with weights as assigned above)
of all the graphs $G[\n]$ and $G'[\n]$:

\begin{theorem}
  \label{thm.Potts_genfn}
Let $G=(V,E)$ be a simple graph,
and let $\bv = (v_e)_{e \in E}$ and $\w = (w_i)_{i \in V}$ be indeterminates.
We then have the following exponential generating functions:
\begin{eqnarray}
   \sum_{\n} Z_{G[\n]}(q, \bv) \, {\x^\n \over \n!}
   & = &
   \left(
   \sum_{\n} \, \Bigl( \prod_{ij \in E}  (1+v_{ij})^{n_i n_j} \Bigr)
             \, {\x^\n \over \n!}
   \right) ^{\! q}
       \label{eq.potts.genfn_1}  \\[4mm]
   \sum_{\n} Z_{G'[\n]}(q, \bv, \w) \, {\x^\n \over \n!}
   & = &
   \left(
   \sum_{\n} \, \Bigl( \prod_{ij \in E}  (1+v_{ij})^{n_i n_j} \Bigr)
             \, \Bigl( \prod_{i \in V}  (1+w_i)^{n_i (n_i-1)/2} \Bigr)
             \, {\x^\n \over \n!}
   \right) ^{\! q}
        \nonumber \\
        \label{eq.potts.genfn_2}
\end{eqnarray}
where all sums run over $\n \in \N^V$.
\end{theorem}

We remark that the expression in large parentheses on the right-hand side
of \reff{eq.potts.genfn_1} [resp.\ \reff{eq.potts.genfn_2}]
is the grand partition function for a lattice gas on the graph $G$
[resp.\ on the graph $G^\circ$ obtained from $G$
 by adjoining a loop at each vertex]
in which arbitrary nonnegative integer occupation numbers
$\n = (n_i)_{i \in V}$ are allowed,
with fugacities $\x$ on the vertices
and two-particle Boltzmann weights $1+v_{ij}$ on edges $ij$
[and two-particle Boltzmann weights $1+w_i$ on the loops].
See e.g.\ \cite{Scott-Sokal_lovasz,Scott-Sokal_lovasz_shortpaper}
for definitions concerning lattice gases.

When $G=K_1$ (the graph with one vertex and no edges),
\reff{eq.potts.genfn_2} reduces to the well-known
\cite{Tutte_67,Biggs_93,Gessel_95,Gessel_96,Welsh_00}
exponential generating function for the
Tutte polynomials of the complete graphs $K_n$:\footnote{
   Unfortunately, this formula is usually written in terms of the
   classical Tutte polynomial $T_G(x,y)$, which is related to 
   $Z_G(q,v)$ by a change of variables \cite[section~2.5]{Sokal_bcc2005}
   that obscures the very different roles played by $q$ and $v$.
}
\begin{equation}
   \sum_{n=0}^\infty Z_n(q,v) \, {x^n \over n!}
   \;=\;
   \left( \sum_{n=0}^\infty (1+v)^{n(n-1)/2} \, {x^n \over n!} \right) ^{\! q}
   \;.
  \label{eq.ZKn.genfn}
\end{equation}
When $G=K_2$, \reff{eq.potts.genfn_1} reduces to the
exponential generating function for the
Tutte polynomials of the complete bipartite graphs $K_{n_1,n_2}$,
\begin{equation}
   \sum_{n_1,n_2=0}^\infty Z_{n_1,n_2}(q,v) \, {x^{n_1} \over {n_1}!}
                                            \, {y^{n_2} \over {n_2}!}
   \;=\;
   \left( \sum_{n_1,n_2=0}^\infty (1+v)^{n_1 n_2} \, {x^{n_1} \over {n_1}!}
                                                  \, {y^{n_2} \over {n_2}!}
   \right) ^{\! q}
   \;,
  \label{eq.ZKn1n2.genfn}
\end{equation}
which in turn specializes for $v=-1$ to the well-known
\cite[Exercise~5.6, pp.~73 and 107--108]{Stanley_99}
exponential generating function for the
chromatic polynomials of the complete bipartite graphs\footnote{
   \reff{eq.PKn1n2.genfn} can alternatively be derived
   from the explicit expression \cite{Swenson_73}
   for $P_{K_{n_1,n_2}}(q)$
   by using a well-known Stirling-number identity \cite[eq.~(7.49)]{Graham_94}
   together with the binomial formula.
},
\begin{equation}
   \sum_{n_1,n_2=0}^\infty P_{K_{n_1,n_2}}(q) \, {x^{n_1} \over {n_1}!}
                                              \, {y^{n_2} \over {n_2}!}
   \;=\;
   (e^x + e^y - 1)^q
   \;.
  \label{eq.PKn1n2.genfn}
\end{equation}
More generally, when $G=K_r$,
\reff{eq.potts.genfn_1} gives the
exponential generating function for the
Tutte polynomials of the complete $r$-partite graphs $K_{n_1,\ldots,n_r}$.

Note also that if we only want to know $Z_{G[\n]}(q, \bv)$ for $\n \le \one$
--- that is, we want to know $Z_{G[W]}(q, \bv)$ for subsets $W \subseteq V$ ---
then it suffices to get correct the terms $\n \le \one$
on the right-hand side of \reff{eq.potts.genfn_1};
all the other terms can be dropped or altered arbitrarily.
It follows that
\begin{equation}
   Z_{G[W]}(q, \bv)   \;=\;  [\x^{\one_W}]
   \left(
   \sum_{\n \le \one} \, \Bigl( \prod_{ij \in E}  (1+v_{ij})^{n_i n_j} \Bigr)
             \, {\x^\n \over \n!}
   \right) ^{\! q}
   \;.
 \label{eq.potts.genfn.nle1.1}
\end{equation}
Note that the expression in large parentheses is now the
grand partition function for a lattice gas with
{\em hard-core self-repulsion}\/,
i.e.\ at most one atom can occupy each site.
In the special case $\bv = -\one$, this becomes
\begin{equation}
   P_{G[W]}(q)   \;=\;  [\x^{\one_W}] \, I_G(\x)^q 
   \;,
 \label{eq.potts.genfn.nle1.2}
\end{equation}
where $I_G(\x)$ is the multivariate generating polynomial
for {\em independent sets}\/ of vertices in $G$.
The formula \reff{eq.potts.genfn.nle1.2} was found previously
by Lass \cite[Proposition~3.1]{Lass_01}.

Since the formula \reff{eq.potts.genfn_1} for $G[\n]$ is just a specialization
of the formula \reff{eq.potts.genfn_2} for $G'[\n]$,
obtained by taking $\w = \zero$, it suffices to prove the latter.

\firstproofof{Theorem~\ref{thm.Potts_genfn}}
Since the coefficient of $\x^\n$ on each side of \reff{eq.potts.genfn_2}
is a polynomial in $q$, it suffices to prove the identity
for infinitely many values of $q$ in the ring $R$;
in particular, it suffices to prove it for all positive integers $q$.
By Theorem~\ref{thm.FK}, for integer $q \ge 1$ we have
\begin{subeqnarray}
   Z_{G'[\n]}(q, \bv, \w)
   & = &
   \sum_{ \sigma \colon\, V(G'[\n]) \to [q]}
   \;\, \prod_{rs \in E(G'[\n])}  \,
      \biggl[ 1 + v_{rs} \delta\bigl(\sigma(r), \sigma(s)\bigr) \biggr]
   \\[4mm]
   & = &
   \!
   \sum_{ \sigma \colon\, V(G'[\n]) \to [q]}
   \; \Biggl(\, \prod_{ij \in E} \prod_{\begin{scarray}
                                        \alpha \in [n_i] \\
                                        \beta \in [n_j]
                                     \end{scarray}}
      \biggl[ 1 + v_{ij} \delta\bigl(\sigma(i,\alpha), \sigma(j,\beta)\bigr)
      \biggr]
      \Biggr)   \;\times
         \nonumber \\
   &  & \qquad\qquad\quad
   \; \Biggl(\, \prod_{i \in V} \prod_{\begin{scarray}
                                        \alpha,\beta \in [n_i] \\
                                        \alpha < \beta
                                    \end{scarray}}
      \biggl[ 1 + w_i \delta\bigl(\sigma(i,\alpha), \sigma(i,\beta)\bigr)
      \biggr]
      \Biggr) \qquad
  \slabel{eq.ZGn}
\end{subeqnarray}
Now, because of the symmetries in the construction of $G'[\n]$,
the summand in \reff{eq.ZGn} depends only on the numbers
$\{n_{i,\tau}\}_{i \in V, \, \tau \in [q]}$ defined by
\begin{equation}
   n_{i,\tau}  \;=\;
   \#\{\alpha \in [n_i] \colon\; \sigma(i,\alpha) = \tau \}
   \;,
\end{equation}
and we have
\begin{eqnarray}
   Z_{G'[\n]}(q, \bv, \w)
   & = &
   \!\!\!\!
   \sum_{\begin{scarray}
            \{n_{i,\tau}\} \\
            \sum_{\tau=1}^q n_{i,\tau} = n_i \, \forall i
         \end{scarray}}
   \left( \prod_{i \in V}  {n_i \choose n_{i,1}, \ldots, n_{i,q}}
   \right)
   \left( \prod_{ij \in E} \, \prod_{\tau=1}^q
          (1+v_{ij})^{n_{i,\tau} n_{j,\tau}}
   \right)
   \;\times
        \nonumber \\
   & &  \qquad\qquad\qquad\;
   \left( \prod_{i \in V} \, \prod_{\tau=1}^q
           (1+w_i)^{n_{i,\tau} (n_{i,\tau}-1)/2}
   \right)
   \;.
\end{eqnarray}
Now multiply by $\x^\n/\n!$ and sum over $\n$;
the sums over the $n_{i,\tau}$ now become unrestricted,
and the sums for the different values of $\tau$ decouple;
we end up with the product of $q$ terms each of which is
\begin{equation}
   \sum_{\n} \, \Bigl( \prod_{ij \in E}  (1+v_{ij})^{n_i n_j} \Bigr)
             \, \Bigl( \prod_{i \in V}  (1+w_i)^{n_i (n_i-1)/2} \Bigr)
             \, {\x^\n \over \n!}
   \;.
\end{equation}
\qed

It is instructive to give an alternate proof of Theorem~\ref{thm.Potts_genfn}
in which $q$ is treated directly as an indeterminate:

\secondproofof{Theorem~\ref{thm.Potts_genfn}}
We begin by proving the special case \reff{eq.ZKn.genfn}.
{}From the binomial formula we have
\begin{subeqnarray}
   \left( 1 \,+\, \sum_{n=1}^\infty (1+v)^{n(n-1)/2} \, {x^n \over n!}
   \right) ^{\! q}
   & = &
   1 \,+\, \sum_{b=1}^\infty  {q^{\underline{b}} \over b!}
      \left( \sum_{n=1}^\infty (1+v)^{n(n-1)/2} \, {x^n \over n!}
      \right) ^{\! b}
         \\[2mm]
   & = &
   1 \,+\, \sum_{b=1}^\infty  {q^{\underline{b}} \over b!}
      \sum_{n_1,\ldots,n_b \ge 1}
          {x^{\sum n_j} \over \prod n_j!} \,
      \prod_{j=1}^b (1+v)^{n_j(n_j-1)/2}
      \;.
   \nonumber \\
  \label{eq.secondproof.0}
\end{subeqnarray}
Now write $n = \sum_{j=1}^b n_j$,
and consider $n_1,\ldots,n_b$ as the block sizes
in an {\em ordered}\/ partition $\pi = (\pi_1,\ldots,\pi_b)$
of $[n]$ into $b$ nonempty blocks.
The number of ordered partitions with those block sizes is precisely the
multinomial coefficient $n! / \prod\limits_{j=1}^b n_j!$,
so \reff{eq.secondproof.0} is equal to
\begin{equation}
   1 \,+\, \sum_{n=1}^\infty  {x^n \over n!}
           \sum_{\pi \in \Pi'_n}  {q^{\underline{|\pi|}} \over |\pi|!}
           \prod_{j=1}^{|\pi|}  (1+v)^{|\pi_j| (|\pi_j|-1)/2}
   \;,
\end{equation}
where $\Pi'_n$ denotes the set of ordered partitions of $[n]$,
and $|\pi|$ denotes the number of blocks in the partition $\pi$.
Passing to {\em unordered}\/ partitions absorbs the factor $1/|\pi|!$
(as these are partitions of a {\em labelled}\/ set),
and we get
\begin{equation}
   1 \,+\, \sum_{n=1}^\infty  {x^n \over n!}
           \sum_{\pi \in \Pi_n}  q^{\underline{|\pi|}}
           \prod_{j=1}^{|\pi|}  (1+v)^{|\pi_j| (|\pi_j|-1)/2}
   \;.
  \label{eq.secondproof.1}
\end{equation}
The factor $(1+v)^{|\pi_j| (|\pi_j|-1)/2}$ corresponds to a sum over
(simple undirected) graphs on the vertex set $\pi_j$,
with a weight $v$ for each edge.
We can now reinterpret \reff{eq.secondproof.1}
as a sum over (simple undirected) graphs $H$ on the vertex set $[n]$,
with a weight $v^{|E(H)|}$, together with a sum over partitions
$\pi \in \Pi_n$ that are {\em compatible}\/ with $H$
in the sense that each block of $\pi$ is a union of (vertex sets of)
connected components of $H$.
If $H$ has $r$ connected components, then these can be grouped
into $b$ nonempty blocks in $\stirlingsubset{r}{b}$ ways
(where $\stirlingsubset{r}{b}$ is a Stirling number of the second kind);
and it is well known \cite[eq.~(6.10)]{Graham_94} that
\begin{equation}
   \sum_{b=1}^r  \stirlingsubset{r}{b} \, q^{\underline{b}}
   \;=\;
   q^r
   \;.
 \label{eq.stirlingidentity}
\end{equation}
Therefore, each graph $H$ gets a weight $q^{k(H)} v^{|E(H)|}$,
and the sum over graphs gives $Z_n(q,v)$.
This proves \reff{eq.ZKn.genfn}.

The more general formulae \reff{eq.potts.genfn_1}/\reff{eq.potts.genfn_2}
are derived by a variant of this proof in which
sums over integers $n \in \N \setminus \{0\}$
are replaced by sums over multi-indices $\n \in \N^V \setminus \{\zero\}$.
In the one-dimensional case, we were considering partitions of $[n]$,
a collection of $n$ distinct objects of a single type.
We are now working with multi-indices in $\N^V$, and so we will be 
partitioning collections of objects of $|V|$ different types.
We begin with some notation.
Given $\n=(n_i)_{i\in V}\in\N^V$, we write $[\n]$ for the set of ordered 
pairs $\{(i,\alpha) \colon\, i\in V \hbox{ and } \alpha\in[n_i]\}$.
Thus $|[\n]|=\sum_{i\in V}n_i$.
We shall think of the subset $[\n]_i:=\{ (i,\alpha) \colon\, \alpha\in[n_i] \}$
as a copy of $[n_i]$ with type $i$;
thus, $[\n]$ is a collection of $|\n|$ objects,
with $n_i$ objects of type $i$ for each $i\in V$.
Note that $G'[\n]$ is simply the complete graph on the vertex set $[\n]$.
Now, given a partition $\bpi=(\bpi_1,\ldots,\bpi_k)$ of $[\n]$
(with all $\bpi_j$ nonempty), we write $|\bpi|=k$.
The {\em block sizes}\/ $\B(\bpi_j)$ are 
elements of $\N^V\setminus\{\mathbf 0\}$ defined by
\begin{equation}
   \B(\bpi_j) \;=\;  (\bpi_j\cap[\n]_i)_{i\in V}
   \;.
\end{equation}
Thus $\B(\bpi_j)$ lists the number of objects of each type contained in 
the $j$th block of $\bpi$.

Now let us define the weight
\begin{equation}
   W(\n) \;:=\;
  \left(\prod_{ii'\in E}(1+v_{ii'})^{n_in_{i'}}\right)
  \left(\prod_{i\in V}(1+w_i)^{n_i(n_i-1)/2}\right)
   \;.
\end{equation}
The right-hand side of \reff{eq.potts.genfn_2} is then
\begin{subeqnarray}
   \left(1+\sum_{\n\ne\mathbf0} W(\n)\frac{\x^{\n}}{\n!}\right)^q
   & = &
   1 \,+\, \sum_{b\ge1}\frac{q^{\underline b}}{b!}
    \left(\sum_{\n\ne\mathbf0}W(\n)\frac{\x^{\n}}{\n!}\right)^b
       \\
   & = &
   1 \,+\, \sum_{b\ge1}\frac{q^{\underline b}}{b!}
     \sum_{\n_1,\ldots,\n_b\ne\mathbf0}
     \frac{\x^{\sum\n_j}}{\prod \n_j!} \prod_{j=1}^b  W(\n_j)  \;.
\label{insertA}
\end{subeqnarray}
Now let $\n=\sum_{j=1}^b\n_j$, and consider $\n_1,\ldots,\n_b$
as the block sizes in 
an ordered partition $\bpi$ of $[\n]$ into $b$ nonempty blocks.
The number of ordered partitions with those block sizes is the
(multidimensional) multinomial coefficient $\n! / \prod\limits_{j=1}^b \n_j!$,
so \reff{insertA} is equal to
\begin{equation}
   1 \,+\, \sum_{\n\ne\mathbf0}\frac{\x^\n}{\n!}
     \sum_\bpi\frac{q^{\underline{|\bpi|}}}{|\bpi|!}
  \prod_{j=1}^{|\bpi|}W(\B(\bpi_j))
   \;,
\end{equation}
where the second sum is over ordered partitions of $[\n]$.
Passing to {\em unordered}\/ partitions, this becomes
\begin{equation}
  \label{insertB}
1 \,+\, \sum_{\n\ne\mathbf0}\frac{\x^\n}{\n!}\sum_\bpi 
q^{\underline{|\bpi|}} \prod_{j=1}^{|\bpi|}W(\B(\bpi_j))  \;.
\end{equation}
The factor $W(\B(\bpi_j))$ corresponds to a sum over graphs with vertex 
set $\bpi_j$, such that edges within $[\n]_i$ have weight $w_i$
and edges between $[\n]_i$ and $[\n]_{i'}$ have weight $v_{ii'}$.
As before, we can reinterpret \reff{insertB} as a sum over 
graphs $H$ with vertex set $[\n]$, with these same edge weights.
The final step of the argument --- summing over partitions $\bpi$
that are compatible with $H$ --- is identical to that in the
one-dimensional case:  it makes no difference whether the vertex set
is $[n]$ or $[\n]$.
We therefore see that each graph $H$ gets weight $q^{k(H)}$ times 
the product of edge weights, and we are done.
\qed

\subsection{Consequences for \mbox{\boldmath$Z_G(q,v)$}}

In view of Theorem~\ref{thm.Potts_genfn}, we can apply all the results
of Section~\ref{sec2} with the identifications
\begin{subeqnarray}
   a_{\n}(q)      & = &  Z_{G'[\n]}(q,\bv,\w)   \\[2mm]
   \ahat_{\n}(q)  & = &  \Zhat_{G'[\n]}(q,\bv,\w)   \\[2mm]
   c_{\n}         & = &  C_{G'[\n]}(\bv,\w)   \\[2mm]
   a_{\n}(1)      & = &
         \left( \prod_{ij \in E}  (1+v_{ij})^{n_i n_j} \right)
         \left( \prod_{i \in V}  (1+w_i)^{n_i (n_i-1)/2} \right)
 \label{eq.complete.ident}
\end{subeqnarray}
We refrain from writing out all the formulae,
which are obtained by substituting (\ref{eq.complete.ident}a--d)
into \reff{eq.partition2.an.cn},
\reff{eq.partition2.inverse.1}/\reff{eq.partition2.inverse.2},
\reff{eq.partition2.rq.q}--\reff{eq.partition2bis.q1hat.q2},
\reff{eq.id1.multidim}, \reff{eq.id2.multidim},
\reff{eq.recursion.an.cn.multidim},
\reff{eq.recursion.anq1.anq2.multidim},
\reff{eq.id1.t.multidim} and \reff{eq.id2.t.multidim}.
Let us simply show a few important results that arise from
the partition formula \reff{eq.partition2.rq.q}/\reff{eq.partition2bis.q1.q2}
and from the convolution formula
\reff{eq.id1.multidim}--\reff{eq.id2.multidim}.

\subsubsection{Partition formulae}  \label{subsubsec.partitions}

By specializing \reff{eq.partition2.rq.q}/\reff{eq.partition2bis.q1.q2}
to $\n = \one$,
we can obtain a partition formula that expresses
the multivariate Tutte polynomial of $G$ at $q=q_2$
in terms of the multivariate Tutte polynomials
of induced subgraphs of $G$ at $q=q_1$,
for {\em arbitrary}\/ choices of $q_1$ and $q_2$:

\begin{proposition}
   \label{prop.ZG.partitions.q1q2}
Let $G=(V,E)$ be a finite loopless graph.  Then
\begin{equation}
   Z_G(q_2, \bv)
   \;=\;
   \sum\limits_{\pi \in \Pi(V)}
      (q_2/q_1)^{\underline{|\pi|}}
      \, \prod_{B\in \pi} \, Z_{G[B]}(q_1,\bv) \;,
 \label{eq.prop.ZG.partitions.q1q2}
\end{equation}
where the sum runs over (unordered) partitions $\pi$ of the vertex set $V$,
the product runs over blocks $B$ of $\pi$,
and $G[B]$ denotes the induced subgraph of $G$ on the vertex set $B$.\footnote{
   Readers who dislike dividing by indeterminates can reinterpret this
   formula by writing $r = q_2/q_1$ and hence $q_2 = rq_1$,
   as in \reff{eq.partition2.rq.q}.
}
Equivalently we can write
\begin{equation}
   Z_G(q_2, \bv)
   \;=\;
   \sum\limits_{\pi \in \Pi(V)}
      \, \Biggl( \, \prod_{j=0}^{|\pi|-1} (q_2 - jq_1) \, \Biggr)
      \, \prod_{B\in \pi} \, \Zhat_{G[B]}(q_1,\bv)
 \label{eq.prop.ZG.partitions.q1q2.alt}
\end{equation}
or
\begin{equation}
   \Zhat_G(q_2, \bv)
   \;=\;
   \sum\limits_{\pi \in \Pi(V)}
      \, \Biggl( \, \prod_{j=1}^{|\pi|-1} (q_2 - jq_1) \, \Biggr)
      \, \prod_{B\in \pi} \, \Zhat_{G[B]}(q_1,\bv) \;.
 \label{eq.prop.ZG.partitions.q1q2.alt2}
\end{equation}
\end{proposition}

\proof
Specializing \reff{eq.partition2.rq.q}
[or equivalently \reff{eq.partition2bis.q1.q2}]
to $\n = \one$,
we obtain a sum over ordered partitions $\pi = (\pi_1,\ldots,\pi_\ell)$
of $V$ into $\ell$ nonempty blocks.
Passing from ordered to unordered partitions,
we get a factor $\ell!$ that cancels the $1/\ell!$
in \reff{eq.partition2.rq.q}/\reff{eq.partition2bis.q1.q2}.
\qed

Proposition~\ref{prop.ZG.partitions.q1q2} is very powerful
because we are free to choose $q_1$ as we please;
then $Z_G(q_2, \bv)$ can be written in terms of $Z_{G[B]}(q_1,\bv)$.
The two most important special cases are $q_1 = 0$ and $q_1 = 1$.\footnote{
   The case $q_1 = -1$ is also of interest, as it is related to
   acyclic orientations.  See Lass \cite{Lass_01} for some interesting
   material that may be relevant in the present context.
}

\medskip

{\bf Case \mbox{\boldmath$q_1=0$}.}
Taking the limit $q_1 \to 0$ in \reff{eq.prop.ZG.partitions.q1q2},
or equivalently just setting $q_1 = 0$
in \reff{eq.prop.ZG.partitions.q1q2.alt},
we obtain
\begin{equation}
   Z_G(q, \bv)
   \;=\;
   \sum\limits_{\pi \in \Pi(V)}
      q^{|\pi|}
      \, \prod_{B\in \pi} \, C_{G[B]}(\bv) \;.
 \label{eq.prop.ZG.partitions.0}
\end{equation}
[This formula can alternatively be derived from \reff{eq.partition2.an.cn}
 by the same reasoning as used in proving
 Proposition~\ref{prop.ZG.partitions.q1q2}.]
Of course, \reff{eq.prop.ZG.partitions.0}
also has an obvious direct combinatorial proof,
based on the definition \reff{def.ZG}:
first we classify subsets $A \subseteq E$
according to the partition $\pi$ of the vertex set $V$
that is induced by the connected components of the subgraph $(V,A)$;
then we sum over ways of connecting up each component.

\medskip

{\bf Case \mbox{\boldmath$q_1=1$}.}
This special case is slightly less obvious,
and it is worth stating it explicitly:

\begin{corollary}
   \label{cor0.ZG.partitions}
Let $G=(V,E)$ be a finite loopless graph.  Then
\begin{eqnarray}
   Z_G(q, \bv)
   & = &
   \sum\limits_{\pi \in \Pi(V)}
      q^{\underline{|\pi|}}
      \, \prod_{B\in \pi} \, \prod_{e \in E(G[B])}  (1+v_e)
 \label{eq.prop.ZG.partitions.1}
\end{eqnarray}
where the sum runs over (unordered) partitions $\pi$ of the vertex set $V$,
the outermost product runs over blocks $B$ of $\pi$,
and the innermost product runs over edges in the induced subgraph $G[B]$.
In particular,
\begin{eqnarray}
   C_G(\bv)
   & = &
   \sum\limits_{\pi \in \Pi(V)}
      (-1)^{|\pi|-1} \: (|\pi| -1)!
      \, \prod_{B\in \pi} \, \prod_{e \in E(G[B])}  (1+v_e)
 \label{eq.prop.ZG.partitions.2}
      \\[2mm]
   Z_G(-1, \bv)
   & = &
   \sum\limits_{\pi \in \Pi(V)}
      (-1)^{|\pi|} \: |\pi|!
      \, \prod_{B\in \pi} \, \prod_{e \in E(G[B])}  (1+v_e)
 \label{eq.prop.ZG.partitions.3}
\end{eqnarray}
where $|\pi|$ denotes the number of blocks in $\pi$.
\end{corollary}

Specializing \reff{eq.prop.ZG.partitions.1}--\reff{eq.prop.ZG.partitions.3}
to $v_e = -1$ for all edges $e$, we have:

\begin{corollary}
   \label{cor.ZG.partitions}
Let $G=(V,E)$ be a finite loopless graph with $|V|=n$.  Then:
\begin{eqnarray}
   P_G(q)  & = &  \sum_{k=1}^n  q^{\underline{k}} \, Q_G(k)
      \label{eq.cor.ZG.partitions.1}  \\[2mm]
   P'_G(0) \;=\; C_G(-1)  & = &  \sum_{k=1}^n  (-1)^{k-1} \, (k-1)! \: Q_G(k)
      \label{eq.cor.ZG.partitions.2}  \\[2mm]
   P_G(-1)  & = &  \sum_{k=1}^n  (-1)^k \, k! \: Q_G(k)
      \label{eq.cor.ZG.partitions.3}
\end{eqnarray}
where $Q_G(k)$ denotes the number of partitions of $V$
into $k$ nonempty {\em independent} subsets.
\end{corollary}

Of course, \reff{eq.cor.ZG.partitions.1} is well known
and has a trivial direct proof by counting proper $q$-colorings;
and \reff{eq.cor.ZG.partitions.2}/\reff{eq.cor.ZG.partitions.3}
are immediate corollaries of \reff{eq.cor.ZG.partitions.1}.


Let us conclude by giving two direct combinatorial proofs of the fundamental
Proposition~\ref{prop.ZG.partitions.q1q2}:
one using the subgraph representation \reff{def.ZG},
and the other using the coloring representation \reff{eq.FK.identity}.

\secondproofof{Proposition~\ref{prop.ZG.partitions.q1q2}}
Use the partition formula \reff{eq.prop.ZG.partitions.0}
[which itself had a simple direct combinatorial proof]
to expand $Z_G(q_2,\bv)$ on the left-hand side of \reff{eq.prop.ZG.partitions.q1q2}
and to expand each $Z_{G[B]}(q_1,\bv)$ on the right-hand side.
We can then interpret the double sum on the right-hand side
as a sum over partitions $\pi$ of $V$, with weight $q_1$ for each block,
together with a sum over ways of grouping these blocks into nonempty groups,
with a weight $(q_2/q_1)^{\underline{k}}$ if we have $k$ groups.
Using the identity
\begin{equation}
   \sum_{k=1}^n  \stirlingsubset{n}{k} \, (q_2/q_1)^{\underline{k}}
   \;=\;
   (q_2/q_1)^n
 \label{eq.stirling2}
\end{equation}
where $n = |\pi|$, we get exactly the expression on the left-hand side.
\qed

It will not escape the reader's notice that the argument using
\reff{eq.stirling2} is identical to the one given earlier
using \reff{eq.stirlingidentity}.

\thirdproofof{Proposition~\ref{prop.ZG.partitions.q1q2}}
Since both sides of \reff{eq.prop.ZG.partitions.q1q2}
are polynomials in $q_1$ and $q_2$, it suffices to prove the identity for
all pairs of positive integers $q_1,q_2$ such that $r = q_2/q_1$
is also an integer.
We shall do this using Theorem~\ref{thm.FK}.
Consider a (not-necessarily-proper) coloring $\sigma \colon\, V \to [q_2]$.
Fix a bijection $[q_2] \simeq [q_1] \times [r]$
and write $\sigma = (\tau,\psi)$ where $\tau \colon\, V \to [q_1]$
and $\psi \colon\, V \to [r]$.
We then trivially have
\begin{subeqnarray}
   Z_G(q_2, \bv)
   & = &
   \sum_{ \sigma \colon\, V \to [q_2]}
   \; \prod_{e=ij \in E}  \,
      \biggl[ 1 + v_e \delta\bigl(\sigma(i), \sigma(j)\bigr) \biggr]
          \\[2mm]
   & = &
   \sum_{ \psi \colon\, V \to [r]}  \;
   \sum_{ \tau \colon\, V \to [q_1]}
   \; \prod_{e=ij \in E}  \,
      \biggl[ 1 + v_e \delta\bigl(\tau(i), \tau(j)\bigr)
                      \delta\bigl(\psi(i), \psi(j)\bigr) \biggr]
   \,. \qquad \quad
\end{subeqnarray}
The map $\psi$ induces a partition $\pi_\psi$ of $V$
into its nonempty color classes,
and the remaining sum over $\tau$ gives (by Theorem~\ref{thm.FK} again)
$\prod_{B \in \pi_\psi} Z_{G[B]}(q_1,\bv)$.
On the other hand, each partition $\pi$ of $V$ arises in this way from
$r^{\underline{|\pi|}}$ different colorings $\psi$,
since we have $r$ choices to color the lexicographically first block of $\pi$,
$r-1$ choices for the lexicographically second block, and so forth.
This proves \reff{eq.prop.ZG.partitions.q1q2}.
\qed

{\bf Important note:}
The formulae in this subsection
--- notably \reff{eq.prop.ZG.partitions.2}
together with the $q=1$ case of \reff{eq.prop.ZG.partitions.0} ---
are strongly reminiscent of M\"obius inversion
on the lattice of partitions $\Pi(V)$.
This is not an accident.
In fact, the formula \reff{eq.prop.ZG.partitions.q1q2.alt2}
is a special case of a ``$q_1$--$q_2$ generalization''
of M\"obius inversion that we present in Appendix~\ref{app.Mobius}.
We are grateful to an anonymous referee for drawing our attention
to the connection of these formulae with M\"obius inversion
and for challenging us to find the underlying general principle.

\subsubsection{Convolution formulae}

In a similar way we can obtain convolution formulae
for the multivariate Tutte polynomial by specializing
\reff{eq.id1.multidim}--\reff{eq.id2.multidim}:


\begin{proposition}
   \label{prop.ZG.recursions}
Let $G=(V,E)$ be a finite loopless graph.
Then
\begin{equation}
   Z_G(q_1+q_2,\bv)  \;=\;
   \sum_{ W \subseteq V }
   Z_{G[W]}(q_1,\bv) \, Z_{G[V \setminus W]}(q_2,\bv)
  \;.
  \label{eq.lin.0}
\end{equation}
Moreover, for each $i \in V$ we have
\begin{equation}
   \Zhat_G(q_1+q_2,\bv)  \;=\;
   \sum_{\begin{scarray}
            W \subseteq V \\
            W \ni i
         \end{scarray}
        }
   \Zhat_{G[W]}(q_1,\bv) \, Z_{G[V \setminus W]}(q_2,\bv)
  \label{eq.lin.1}
\end{equation}
and in particular
\begin{equation}
   Z_G(q,\bv)  \;=\;  \sum_{\begin{scarray}
                               W \subseteq V \\
                               W \ni i
                            \end{scarray}
                           }
   q C_{G[W]}(\bv) \, Z_{G[V \setminus W]}(q,\bv)
   \;.
   \label{eq.lin.2}
\end{equation}
We also have
\begin{equation}
   |V| \, Z_G(q_2,\bv)
   \;=\,
   \sum_{\emptyset \neq W \subseteq V}
      \Bigl[ (q_1+q_2)|W| \,-\, q_1 |V| \Bigr] \,
   \Zhat_{G[W]}(q_1,\bv) \, Z_{G[V \setminus W]}(q_2,\bv)
   \;.  \qquad
   \label{eq.lin.2a}
\end{equation}
\end{proposition}

\proof
\reff{eq.lin.0} and \reff{eq.lin.1} are just
\reff{eq.id1.multidim} and \reff{eq.id2.multidim}, respectively,
specialized to $\n = \one$.
Further specializing the latter to $q_1 = 0$, we obtain \reff{eq.lin.2}.
Finally, \reff{eq.lin.2a} is obtained by subtracting
$(q_1+q_2) \times \hbox{\reff{eq.lin.1}}$ from \reff{eq.lin.0}
and summing over $i \in V$.
\qed

For completeness, let us give a simple direct proof of \reff{eq.lin.1}.
[\reff{eq.lin.0} is similar but easier, and is left as an exercise
 for the reader.]

\altproofof{Proposition~\ref{prop.ZG.recursions}}
{}From the definition \reff{def.ZhatG} we have
\begin{equation}
   \Zhat_G(q_1+q_2, \bv)   \;=\;
   \sum_{A \subseteq E}  (q_1+q_2)^{k(A)-1}  \prod_{e \in A}  v_e
   \;.
\end{equation}
The factor $(q_1+q_2)^{k(A)-1}$ can be handled by choosing,
for each connected component of $(V,A)$
{\em other than}\/ the component containing the distinguished vertex $i$,
to color it either ``1'' or ``2'', with a corresponding factor $q_1$ or $q_2$,
and summing over all such choices.
Now define $W \subseteq V$ to be the vertex set corresponding to
the union of the components colored ``1''
{\em together with}\/ the component containing $i$.
By construction, $W$ is compatible with $A$
in the sense that $A$ has no edges connecting $W$ to $V \setminus W$.
Moreover, every $W \ni i$ compatible with $A$ is obtained exactly once
by some choice of colors ``1'' and ``2''.
On the other hand, the right-hand side of \reff{eq.lin.1}
is given precisely by such a sum over compatible pairs $(W,A)$
[the $\widehat{\hphantom{Z}}$ on $\Zhat_{G[W]}(q_1,\bv)$
 accounts for the fact that the component containing $i$
 carries no color and hence no factor $q_1$].
\qed

We do not know whether the Abel-type extensions
\reff{eq.id1.t.multidim}/\reff{eq.id2.t.multidim},
when applied to the multivariate Tutte polynomial, are of any interest.

\subsection{Specialization to the complete graphs}

Let us now take $G$ to be the complete graph $K_n$,
with all the edge weights $v_e$ set equal to the same value $v$.
We write $Z_n(q,v)$, $\Zhat_n(q,v)$ and $C_n(v)$
for the corresponding polynomials.

\subsubsection{Partition formulae}

Specializing \reff{eq.prop.ZG.partitions.q1q2},
we obtain the general identity
\begin{equation}
   Z_n(q_2,v)
   \;=\;  
   \sum\limits_{\pi \in \Pi_n}  (q_2/q_1)^{\underline{|\pi|}}
           \prod_{i=1}^{|\pi|}  Z_{|\pi_i|}(q_1,v)
   \;,
 \label{eq.Znqv.partitions.q1q2}
\end{equation}
where the sum runs over (unordered) partitions $\pi$ of $[n]$,
say $\pi = \{\pi_1,\ldots,\pi_\ell\}$,
and $|\pi|$ denotes the number of blocks in $\pi$ (i.e., $|\pi| = \ell$).
Special cases of \reff{eq.Znqv.partitions.q1q2}
arise at particular values of $q_1$:

\medskip

{\bf Case \mbox{\boldmath$q_1=0$}.}
Specializing \reff{eq.prop.ZG.partitions.0}, we obtain the trivial identity
\begin{equation}
   Z_n(q,v)
   \;=\;
   \sum\limits_{\pi \in \Pi_n}  q^{|\pi|} \prod_{i=1}^{|\pi|} C_{|\pi_i|}(v)
   \;.
 \label{eq.Znqv.partitions.0}
\end{equation}

\medskip

{\bf Case \mbox{\boldmath$q_1=1$}.}
Specializing \reff{eq.prop.ZG.partitions.1}--\reff{eq.prop.ZG.partitions.3},
we obtain the less trivial relations
\begin{eqnarray}
   Z_n(q, v)
   & = &
   \sum\limits_{\pi \in \Pi_n}
      q^{\underline{|\pi|}}
      \: (1+v)^{n(n-1)/2 - \|\pi\|}
      \label{eq.Znqv.partitions}   \\[2mm]
   C_n(v)
   & = &
   \sum\limits_{\pi \in \Pi_n}
      (-1)^{|\pi|-1} \: (|\pi| -1)!
      \: (1+v)^{n(n-1)/2 - \|\pi\|}
   \label{eq.Cnv.partitions}  \\[2mm]
   Z_n(-1, v)
   & = &
   \sum\limits_{\pi \in \Pi_n}
      (-1)^{|\pi|} \: |\pi|!
      \: (1+v)^{n(n-1)/2 - \|\pi\|}
   \label{eq.Znminusonev.partitions}
\end{eqnarray}
where
\begin{equation}
   \| \pi \|  \;=\;  \sum_{1 \le i<j \le \ell} |\pi_i|  \, |\pi_j|
        \;=\; {1 \over 2} \left( n^2 \,-\, \sum_{i=1}^\ell |\pi_i|^2 \right)
 \label{def.crossedges}
\end{equation}
denotes the number of ``cross-edges'' in $\pi$,
i.e.\ the number of edges in the complete $\ell$-partite graph
with vertex classes $\pi_1,\ldots,\pi_\ell$.
Equivalently,
\begin{equation}
   {n(n-1) \over 2} \,-\, \| \pi \|
   \;=\;
   \sum_{i=1}^\ell {|\pi_i| (|\pi_i| - 1)  \over 2}
 \label{eq.crossedges.2}
\end{equation}
is the number of ``internal edges'' in $\pi$.
The formulae \reff{eq.Znqv.partitions}--\reff{eq.Znminusonev.partitions}
will play an important role in our study of the large-$n$ asymptotics
of $C_n(v)$ and $Z_n(q,v)$ at complex $v$ \cite{Scott-Sokal_cn_bounds}.

\subsubsection{Convolution formulae}

Specializing \reff{eq.lin.0}--\reff{eq.lin.2a}
and exploiting the symmetries of the complete graph, we obtain
\begin{eqnarray}
   Z_n(q_1+q_2,v)
   & = &
   \sum_{k=0}^n {n \choose k} \, Z_k(q_1,v) \, Z_{n-k}(q_2,v)
       \label{eq.lin.0.Kn}  \\[2mm]
   \Zhat_n(q_1+q_2,v)
   & = &
   \sum_{k=1}^n {n-1 \choose k-1} \, \Zhat_k(q_1,v) \, Z_{n-k}(q_2,v)
   \qquad\hbox{for } n \ge 1
       \label{eq.lin.1.Kn}  \\[2mm]
   Z_n(q,v)
   & = &
   \sum_{k=1}^n {n-1 \choose k-1} \, q C_k(v) \, Z_{n-k}(q,v)
   \qquad\hbox{for } n \ge 1
       \label{eq.lin.2.Kn}  \\[2mm]
   Z_n(q_2,v)
   & = &
   \sum_{k=1}^n \left[ {n-1 \choose k-1} (q_1+q_2) \,-\, {n \choose k} q_1
                \right]  \Zhat_k(q_1,v) \, Z_{n-k}(q_2,v)
   \quad\hbox{for } n \ge 1
       \label{eq.lin.2a.Kn}  \nonumber \\
\end{eqnarray}
On the other hand, these are just the one-dimensional identities
\reff{eq.id1}, \reff{eq.id2}, \reff{eq.recursion.an.cn}
and \reff{eq.recursion.anq1.anq2}
applied to the family $a_n(q) = Z_n(q,v)$,
which is indeed of the required form \reff{eq1.2}
by virtue of the generating-function formula \reff{eq.ZKn.genfn}.

The formula \reff{eq.lin.2.Kn} can be used to compute the $Z_n(q,v)$
inductively given the $C_k(v)$.
On the other hand, if we specialize \reff{eq.lin.2.Kn} to $q=1$
and use \reff{ZG.q=1}, we obtain the identity
\begin{equation}
   (1+v)^{n(n-1)/2}  \;=\;
     \sum_{k=1}^{n} {n-1 \choose k-1} \, C_k(v) \, (1+v)^{(n-k)(n-k-1)/2}
   \qquad\hbox{for } n \ge 1
  \label{eq.lin.3.Kn}
\end{equation}
or equivalently
\begin{equation}
   C_n(v)  \;=\;  (1+v)^{n(n-1)/2}  \,-\,
     \sum_{k=1}^{n-1} {n-1 \choose k-1} \, C_k(v) \, (1+v)^{(n-k)(n-k-1)/2}
   \qquad\hbox{for } n \ge 1
   \;,
  \label{eq.lin.3bis.Kn}
\end{equation}
which can be used to compute the $C_n(v)$ inductively {\em ab initio}\/.
We call \reff{eq.lin.3.Kn}/\reff{eq.lin.3bis.Kn}
the ``linear'' identity for $C_n(v)$.
It goes back at least to Leroux \cite[eq.~(3.3)]{Leroux_88}
and is probably much older;
it can be proven by an easy direct argument based on considering
the size $k$ of the connected component containing a fixed vertex.
[The same argument proves \reff{eq.lin.2.Kn}.]

Alternatively, if we specialize \reff{eq.lin.2a.Kn} to $q_1=1$
and use \reff{ZG.q=1}, we obtain
\begin{equation}
   Z_n(q,v)   \;=\;
   \sum_{k=1}^{n}  \left[ {n-1 \choose k-1} (1+q) \,-\, {n \choose k}
                   \right]
     (1+v)^{k(k-1)/2} \, Z_{n-k}(q,v)
   \qquad\hbox{for } n \ge 1
   \;,
 \label{eq.lin.4.Kn}
\end{equation}
which allows us to go directly from
from $Z_n(1,v) = (1+v)^{n(n-1)/2}$ to $Z_n(q,v)$
without passing through $C_n(v)$.
We can also rewrite \reff{eq.lin.4.Kn} as
\begin{eqnarray}
   \Zhat_n(q,v)
   & = &
   (1+v)^{n(n-1)/2}  \,+\,
   \sum_{k=1}^{n-1}  \left[ {n-1 \choose k} q \,- {n-1 \choose k-1}
                     \right]
     (1+v)^{(n-k)(n-k-1)/2} \, \Zhat_{k}(q,v)
        \nonumber \\
   & &
   \hspace*{3.5in}
   \qquad\hbox{for } n \ge 1
   \;,
 \label{eq.lin.4bis.Kn}
\end{eqnarray}
which manifestly generalizes \reff{eq.lin.3bis.Kn}
and reduces to it when $q=0$.

\subsection{A generalization}

Let ${\bf a} = \{a_n\}_{n=0}^\infty$ be an arbitrary sequence
of coefficients belonging to the ring $R$, satisfying $a_0 = 1$,
and define the family of polynomials
\begin{equation}
   Z_n(q; {\bf a})
   \;=\;
   \sum\limits_{\pi \in \Pi_n}  q^{\underline{|\pi|}}
       \prod_{i=1}^{|\pi|} a_{|\pi_i|}
 \label{def.Znqa}
\end{equation}
with the convention $Z_0 = 1$.
Note that $Z_n(1; {\bf a}) = a_n$.
This definition generalizes \reff{eq.Znqv.partitions}/\reff{eq.crossedges.2},
and if we specialize to $a_n = (1+v)^{n(n-1)/2}$ we obtain $Z_n(q,v)$.

A simple counting argument shows that
\begin{equation}
   Z_n(q; {\bf a})
   \;=\;
   \sum_{b=0}^n  {q^{\underline{b}} \over b!}
      \sum_{\begin{scarray}
               n_1, \ldots, n_b \ge 1  \\
               \sum n_i = n
            \end{scarray}}
      \!\!
      {n \choose n_1, \ldots, n_b}  \prod_{i=1}^b a_{n_i}
   \;.
 \label{eq.Znqa.bis}
\end{equation}
Using this to compute the exponential generating function
of the $\{ Z_n(q; {\bf a}) \}$, we find,
after a short calculation using the binomial series, that
%
\begin{equation}
   \sum_{n=0}^\infty {x^n \over n!} \, Z_n(q; {\bf a})
   \;=\;
   \left( \sum_{n=0}^\infty {x^n \over n!} \, a_n \right) ^{\! q}
    \;.
\end{equation}
In other words, the family $\{ Z_n(q; {\bf a}) \}$
is of the form $\{a_n(q)\}$ defined in \reff{eq1.2},
with $a_n(1) = a_n$.
But this should hardly be surprising, as \reff{eq.Znqa.bis}
is simply the one-dimensional case of \reff{eq.partition2bis.q1.q2}
with $q_2 = q$ and $q_1 = 1$.
We have thus come full circle.

Likewise, the recursion
\begin{equation}
   Z_n(q; {\bf a})
   \;=\;
   q \sum_{k=1}^n {n-1 \choose k-1} \, a_k \, Z_{n-k}(q-1; {\bf a})
\end{equation}
can be proven from \reff{def.Znqa}
by an easy direct argument based on considering
the size $k$ of the block of $\pi$ containing some fixed element of $[n]$.
But this identity is nothing other than \reff{eq.id2.bis}
specialized to $q_1 = 1$ and $q_2 = q-1$.

\subsection{A related problem}

Suppose we try a definition analogous to \reff{def.Znqa},
but with ordinary powers $q^{|\pi|}$
in place of falling factorials $q^{\underline{|\pi|}}$.
That is, let ${\bf c} = \{c_n\}_{n=1}^\infty$ be an arbitrary sequence
of coefficients belonging to the ring $R$,
and define the family of polynomials
\begin{equation}
   Y_n(q; {\bf c})
   \;=\;
   \sum\limits_{\pi \in \Pi_n}  q^{|\pi|} \prod_{i=1}^{|\pi|} c_{|\pi_i|}
 \label{def.Yn}
\end{equation}
with the convention $Y_0 = 1$.
Of course, the factor $q^{|\pi|}$ is superfluous,
because we can simply multiply each $c_n$ by $q$,
but it is convenient to keep it explicit.
This definition generalizes \reff{eq.Znqv.partitions.0},
and if we specialize to $c_n = C_n(v)$ we obtain $Z_n(q,v)$.

A simple counting argument shows that
\begin{equation}
   Y_n(q; {\bf c})
   \;=\;
   \sum_{b=0}^n  {q^b \over b!}
      \sum_{\begin{scarray}
               n_1, \ldots, n_b \ge 1  \\
               \sum n_i = n
            \end{scarray}}
      \!\!
      {n \choose n_1, \ldots, n_b}  \prod_{i=1}^b c_{n_i}
   \;.
 \label{eq.Yn.bis}
\end{equation}
Using this to compute the exponential generating function
of the $\{ Y_n(q; {\bf c}) \}$, we find,
after a short calculation using the exponential series, that
%
\begin{equation}
   \sum_{n=0}^\infty {x^n \over n!} \, Y_n(q; {\bf c})
   \;=\;
   \exp\! \left( q \sum_{n=1}^\infty {x^n \over n!} \, c_n \right)
    \;.
\end{equation}
In other words, the family $\{ Y_n(q; {\bf c}) \}$
is again of the form $\{a_n(q)\}$ defined in \reff{eq1.2}.
But this should hardly be surprising, as \reff{eq.Yn.bis}
is simply the one-dimensional case of \reff{eq.partition2.an.cn}
[i.e.\ \reff{eq.partition.an.cn}.]
We have again come full circle.

Likewise, the recursion
\begin{equation}
   Y_n(q; {\bf c})
   \;=\;
   q \sum_{k=1}^n {n-1 \choose k-1} \, c_k \, Y_{n-k}(q; {\bf c})
   \;,
\end{equation}
which can be proven from \reff{def.Yn}
by considering the size $k$ of the block of $\pi$
containing some fixed element of $[n]$,
is nothing other than \reff{eq.recursion.an.cn}.

\section{Nonlinear identity for the multivariate Tutte polynomial}

%

\subsection{General identity}

In the preceding section we proved the identities
\reff{eq.lin.1} and \reff{eq.lin.2} for the multivariate Tutte polynomial
$Z_G(q,\bv)$, which are based on choosing a single distinguished vertex $i$.
We shall now prove a different identity that is based on choosing
a {\em pair}\/ of distinguished vertices $i,j$.
Unlike the identities discussed in the preceding section,
this one relies on the graphical structure of $Z_G(q,\bv)$
and does {\em not}\/ appear to generalize to arbitrary families
$\{a_n(q)\}$ of the type \reff{eq1.2}.

So let $G = (V,E)$ be a finite undirected graph with $|V| \ge 2$,
and let $i,j \in V$ with $i \neq j$.
We then have the following identity:

\begin{theorem}
   \label{thm.nonlin}
Let $G = (V,E)$ be a finite graph, and let $i,j \in V$ with $i \neq j$.
Then
\begin{equation}
   Z_G(q,\bv)  \;=\;
   \sum_{\begin{scarray}
              W \subseteq V \\
              W \ni i, \, W \not\ni j
         \end{scarray}
        }
   \left[ q - 1 + \prod_{e \in E(W,j)} (1+v_e) \right]
   C_{G[W]}(\bv) \, Z_{G[V \setminus W]}(q,\bv)
 \label{eq.nonlin}
\end{equation}
where $E(W,j)$ denotes the set of all edges with one endpoint in $W$
and the other endpoint at $j$.
In particular,
\begin{equation}
   C_G(\bv)  \;=\;
   \sum_{\begin{scarray}
              W \subseteq V \\
              W \ni i, \, W \not\ni j
         \end{scarray}
        }
   \left[ \prod_{e \in E(W,j)} (1+v_e) \,-\, 1 \right]
   C_{G[W]}(\bv) \, C_{G[V \setminus W]}(\bv)
   \;.
 \label{eq.nonlin.2}
\end{equation}
\end{theorem}

We do not know whether there exists a generalization of \reff{eq.nonlin}
involving $q_1$ and $q_2$ (rather than just $q_1 = 0$).

The proof of Theorem~\ref{thm.nonlin} is based on
looking at the connected component of $i$ in the induced subgraph
where $j$ is deleted:

\proof
Start from the definition
$Z_G(q, \bv) = \sum\limits_{A \subseteq E}
                  q^{k(A)} \prod\limits_{e \in A}  v_e$.
Let $G' = (V,A)$ and $G'' = G' \setminus j$,
and let $W$ be the vertex set of the connected component of $G''$
containing $i$.  Let us now sum over all $A \subseteq E$
that give rise in this way to a specified set $W \subseteq V$,
and let us split this sum into two parts according as
$A$ does or does not contain at least one edge from $j$ to $W$.
If $A$ does not contain such an edge,
then $W$ is a connected component of $G'$,
giving rise to a factor $q$, and we get
$q C_{G[W]}(\bv) \, Z_{G[V \setminus W]}(q,\bv)$.
On the other hand, if $A$ contains at least one such edge,
then $W$ forms part of the connected component of $j$ ($\in V \setminus W$)
in $G'$;  hence there is no factor $q$, but there is a factor $v_e$
for each edge $e$ in the (nonempty) subset $A \cap E(W,j)$.
Summing over all such nonempty subsets of $E(W,j)$, we get
\begin{equation}
   \left[ \prod_{e \in E(W,j)} (1+v_e) \,-\, 1 \right]
   C_{G[W]}(\bv) \, Z_{G[V \setminus W]}(q,\bv)
   \;.
\end{equation}
Putting everything together gives \reff{eq.nonlin}.
Specialization to $q=0$ yields \reff{eq.nonlin.2}.
\qed

There is also a variant of \reff{eq.nonlin} in which
the right-hand side is completely expanded out.
For notational simplicity we shall assume that
$G$ has no loop at the distinguished vertex $j$.
(In the general case one must multiply by a factor $1+v_e$
 for each loop $e$ at $j$.)

\begin{proposition}
   \label{prop.genborgs}
Let $G = (V,E)$ be a finite graph, let $j \in V$,
and suppose that $G$ has no loop at $j$.
Then
\begin{equation}
   \Zhat_G(q,\bv)  \;=\;
   \sum_{\pi \in \Pi(V \setminus \{j\})}
   \: \prod_{B \in \pi}
   \left[ q - 1 + \prod_{e \in E(B,j)} (1+v_e) \right]
   C_{G[B]}(\bv)
 \label{eq.genborgs.1}
\end{equation}
where the sum runs over (unordered) partitions $\pi$ of $V \setminus \{j\}$,
and the product runs over blocks $B$ of $\pi$.
In particular,
\begin{equation}
   C_G(\bv)  \;=\;
   \sum_{\pi \in \Pi(V \setminus \{j\})}
   \: \prod_{B \in \pi}
   \left[ \prod_{e \in E(B,j)} (1+v_e) \,-\, 1 \right]
   C_{G[B]}(\bv)
 \label{eq.genborgs.2}
\end{equation}
\end{proposition}

\proof
In \reff{def.ZhatG}, consider a term $A \subseteq E$,
and let $G' = (V,A)$ and $G'' = G' \setminus j$.
Let $\pi$ be the partition of $V \setminus \{j\}$
into vertex sets of connected components of $G''$.
We can recover $G'$ from $G''$ by adjoining,
for each block $B$ of $\pi$,
zero or more edges from the set $E(B,j)$.
If we adjoin zero edges, we get an extra factor $q$
because $B$ becomes the vertex set of a connected component of $G'$
that is distinct from the component containing $j$;
if we adjoin one or more edges, we get no such factor.
This proves \reff{eq.genborgs.1}.
Specialization to $q=0$ yields \reff{eq.genborgs.2}.
\qed

When $v_e=v$ for all edges $e$,
the formulae \reff{eq.genborgs.1} and \reff{eq.genborgs.2}
can be found in Gessel \cite[Theorems~13 and 10]{Gessel_95}.
The special case of \reff{eq.genborgs.2} in which $v_e=-1$
for all $e$ was also found by Borgs \cite{Borgs_00}
and used recently by him \cite[Lemma~3.2]{Borgs_06}
to bound the complex zeros of chromatic polynomials
(a variant of the proof in \cite{Sokal_chromatic_bounds}).
Indeed, we were led to formulate \reff{eq.genborgs.1}/\reff{eq.genborgs.2}
by meditating on Borgs' special case,
oblivious to the fact that they had already been essentially found by Gessel!

\bigskip

Finally, let us use \reff{eq.nonlin.2} to prove an interesting inequality
concerning the polynomials $C_G(\bv)$.
Let
\begin{equation}
   c(A)  \;=\;  |A| \,-\, |V| \,+\, k(A)
\end{equation}
be the cyclomatic number of the subgraph $(V,A)$,
and let us define the generalized connected sum
\begin{subeqnarray}
   C_G(\bv,\lambda)
   & = &
   \sum_{\begin{scarray}
             A \subseteq E \\
             k(A) = 1
         \end{scarray}}
   \lambda^{c(A)} \prod_{e \in A} v_e
        \\[2mm]
   & = &
   \lambda^{-(|V|-1)} \, C_G(\lambda\bv)
   \;.
        \slabel{eq.CGlambda}
\end{subeqnarray}
Of course, \reff{eq.CGlambda} shows that $C_G(\bv,\lambda)$ contains
no more information than $C_G(\bv)$;
it is simply a convenient way of scaling all the variables $v_e$ simultaneously
while removing a factor $\lambda^{|V|-1}$.
In particular, $C_G(\bv,\lambda)$ interpolates between
the spanning-tree sum ($\lambda=0$) and the connected-spanning-subgraph sum
($\lambda=1$).  We then have the following result
\cite[Remark~2 in Section~4.1]{Sokal_chromatic_bounds}
\cite[Proposition~2.5]{Scott-Sokal_lovasz}:

\begin{proposition}
   \label{prop.ineq.CGlambda}
Let $G=(V,E)$ be a finite graph ($V \neq \emptyset$)
equipped with real edge weights $\bv = (v_e)_{e \in E}$
satisfying $-1 \le v_e \le 0$ for all $e \in E$.
Then
\begin{equation}
   (-1)^{\ell+|V|-1} \, {\partial^\ell \over \partial\lambda^\ell} \,
       C_G(\bv,\lambda)
   \;\ge\;
   0
 \label{eq.prop.ineq.CGlambda}
\end{equation}
on $0 \le \lambda \le 1$, for all integers $\ell \ge 0$.
\end{proposition} 

This inequality was proven in \cite{Sokal_chromatic_bounds,Scott-Sokal_lovasz}
using a ``partitionability'' method going back to
Penrose \cite{Penrose_67}.\footnote{
   For the partitionability method
   (which applies to matroids as well as graphs),
   see also \cite{Bjorner_92,Gessel_96} and the other references mentioned in
   \cite[Section~2.2]{Scott-Sokal_lovasz}.
}
Here we prove it using the recursion \reff{eq.nonlin.2}:

\proofof{Proposition~\ref{prop.ineq.CGlambda}}
By induction on $|V|$.
When $|V|=1$, $G$ necessarily consists of zero or more loops
attached to the sole vertex, so
\begin{equation}
   C_G(\bv,\lambda)  \;=\;  \prod_{e \in E} (1+\lambda v_e)
   \;.
\end{equation}
Using $v_e \le 0$ and $1+\lambda v_e \ge 0$,
it is easy to see that \reff{eq.prop.ineq.CGlambda} holds.

Now assume that $|V| \ge 2$.
Replacing $\bv$ by $\lambda\bv$ in \reff{eq.nonlin.2}
and using \reff{eq.CGlambda}, we obtain
\begin{subeqnarray}
   C_G(\bv,\lambda)
   & = &
   \!\!
   \sum_{\begin{scarray}
              W \subseteq V \\
              W \ni i, \, W \not\ni j
         \end{scarray}
        }
   \!\!
   \lambda^{-1} \! \left[ \prod_{e \in E(W,j)} (1+v_e) \,-\, 1 \right]
   C_{G[W]}(\bv,\lambda) \, C_{G[V \setminus W]}(\bv,\lambda)
         \nonumber \\ \\
   & = &
   \!\!
   \sum_{\begin{scarray}
              W \subseteq V \\
              W \ni i, \, W \not\ni j
         \end{scarray}
        }
   \!\!
   \left( \! \sum_{\emptyset \neq B \subseteq E(W,j)} \!
          \lambda^{|B|-1}  \prod_{e \in B} v_e
   \right)
   C_{G[W]}(\bv,\lambda) \, C_{G[V \setminus W]}(\bv,\lambda)
   \;.
         \nonumber \\
\end{subeqnarray}
Now apply $\partial^\ell/\partial\lambda^\ell$ to both sides:
on the right-hand side we will have terms in which
$\ell_1$~derivatives act on $\sum \lambda^{|B|-1} \prod v_e$,
$\ell_2$ act on $C_{G[W]}(\bv,\lambda)$, and
$\ell_3$ act on $C_{G[V \setminus W]}(\bv,\lambda)$,
where $\ell_1 + \ell_2 + \ell_3 = \ell$.
The induction hypothesis is applicable because $1 \le |W| \le |V|-1$.
So it suffices to check that
\begin{equation}
   (-1)^{\ell_1 -1} \, {\partial^{\ell_1} \over \partial\lambda^{\ell_1}}
   \left( \! \sum_{\emptyset \neq B \subseteq E(W,j)} \!
          \lambda^{|B|-1}  \prod_{e \in B} v_e
   \right)
   \;\ge\;
   0  \;.
 \label{ineq.K2m}
\end{equation}
To prove this, let us order the elements of $E(W,j)$,
i.e.\ consider $E(W,j) \simeq [m]$ where $m = |E(W,j)|$.
Then, by singling out the smallest element in each set $B$, we can write
\begin{equation}
   \sum_{\emptyset \neq B \subseteq [m]} \!
          \lambda^{|B|-1}  \prod_{e \in B} v_e
   \;=\;
   \sum_{i=1}^m  v_i \prod_{j=i+1}^m (1+ \lambda v_j)
   \;.
 \label{eq.K2m.bis}
\end{equation}
We have $v_i \le 0$;
each $\partial/\partial\lambda$ brings down some factor $v_j \le 0$;
and the undifferentiated factors satisfy $1+\lambda v_j \ge 0$.
So \reff{ineq.K2m} holds.
\qed

It is worth remarking that \reff{ineq.K2m} is nothing other than
Proposition~\ref{prop.ineq.CGlambda} specialized to the graph $K_2^{(m)}$
consisting of two vertices connected by $m$ parallel edges.
Moreover, the identity \reff{eq.K2m.bis} is nothing other than
(one version of) the partitionability identity for $K_2^{(m)}$.

\subsection{Specialization to the complete graphs}

Specializing \reff{eq.nonlin} to the complete graph $K_n$
with equal weights $v$, we obtain
\begin{equation}
   Z_n(q,v)  \;=\;
   \sum_{k=1}^{n-1} {n-2 \choose k-1} \, [q+(1+v)^k -1] \,
       C_k(v) \, Z_{n-k}(q,v)
   \qquad\hbox{for } n \ge 2
 \label{eq.nonlin.Kn}
\end{equation}
(see \cite[eq.~(5)]{Gessel_95} for an equivalent formula).
Specializing this to $q=0$, we obtain
\begin{equation}
   C_n(v)  \;=\;
   \sum_{k=1}^{n-1} {n-2 \choose k-1} \, [(1+v)^k - 1] \,
       C_k(v) \, C_{n-k}(v)
   \qquad\hbox{for } n \ge 2
   \;,
 \label{eq.nonlin.Kn.q=0}
\end{equation}
a result apparently first proven by Leroux \cite[eq.~(3.5)]{Leroux_88}
(see also \cite[eq.~(2)]{Gessel_95}).\footnote{
   See also \cite[p.~306, Exercise~4.2.2]{Bergeron_98}
   for a proof of \reff{eq.nonlin.Kn.q=0}
   in species language.
}
And specializing this latter formula to $v=1$, we obtain the identity
of Riordan, Nijenhuis, Wilf and Kreweras
\cite{Riordan_unpub,Nijenhuis_79,Kreweras_80}
for counting connected graphs.
We call \reff{eq.nonlin.Kn.q=0} the ``nonlinear'' identity for $C_n(v)$.

On the other hand, specializing \reff{eq.nonlin.Kn} to $q=1$,
dividing both sides by $(1+v)^{n-1}$ and relabelling $n \to n+1$,
we recover the ``linear'' identity \reff{eq.lin.3.Kn}.

The formula \reff{eq.nonlin.Kn.q=0} has also been derived
in a different context.
Let $T$ be a tree with vertex set $[n]$, rooted at the vertex 1.
An {\em inversion}\/ of $T$ is an ordered pair $(j,k)$ of vertices
such that $j > k > 1$ and $k$ is a descendant of $j$
(i.e., the path from 1 to $k$ passes through $j$).
We define the {\em inversion enumerator for trees}\/
\cite{Mallows_68,Gessel_79,Kreweras_80}
(see also \cite{Gessel_95,Gessel_96,Chauve_00})
to be the polynomial
\begin{equation}
   I_n(y)  \;=\; \sum_{\hbox{\scriptsize trees $T$ on $[n]$}}
                 y^{{\rm inv}(T)}
 \label{def.In}
\end{equation}
where ${\rm inv}(T)$ denotes the number of inversions in $T$.
This polynomial turns out to be related to $C_n(v)$ by the beautiful formula
\cite{Mallows_68,Gessel_79,Kreweras_80,Beissinger_82}
\begin{equation}
   C_n(v)  \;=\;  v^{n-1} I_n(1+v)
   \;.
 \label{eq.Cn.In}
\end{equation}
Now, Mallows, Riordan and Kreweras \cite{Mallows_68,Kreweras_80}
show that $I_n(y)$ satisfies the recursion
\begin{equation}
   I_n(y)  \;=\;  
   \sum_{k=1}^{n-1} {n-2 \choose k-1} \left( \sum_{j=0}^{k-1} y^j \right)
       I_k(y) \, I_{n-k}(y)
   \qquad\hbox{for } n \ge 2
   \;.
 \label{eq.nonlin.In}
\end{equation}
But using \reff{eq.Cn.In}, it is easily seen that \reff{eq.nonlin.In}
is equivalent to \reff{eq.nonlin.Kn.q=0}.

Let us remark, finally, that even without using \reff{def.In},
it follows immediately from the recursion \reff{eq.nonlin.In}
[together with the initial condition $I_1(y) = 1$]
that the polynomials $I_n(y) \equiv C_n(y-1)/(y-1)^{n-1}$
have nonnegative (indeed, strictly positive) integer coefficients.
On the other hand, the nonnegativity of the derivatives of $I_n(y)$
at $y=0$ is also a special case of Proposition~\ref{prop.ineq.CGlambda}.


\appendix
\section{A generalization of M\"obius inversion on the partition lattice}
   \label{app.Mobius}

In this appendix we present an apparently new generalization
of M\"obius inversion on the lattice of partitions of a finite set,
which is inspired by the formulae in Section~\ref{subsubsec.partitions}
and in particular by \reff{eq.prop.ZG.partitions.q1q2.alt2},
and more generally by \reff{eq.partition2bis.q1hat.q2}.
We refer to \cite[Chapter~3]{Stanley_86} for basic facts about posets
and M\"obius inversion.

To begin with, let $P$ be a finite poset.
Then the {\em zeta function}\/ $\zeta$ on $P$ is the function
$\zeta \colon\, P \times P \to \Z$ defined by
\be
   \zeta(x,y)  \;=\;
   \begin{cases}
      1  & \text{if $x \le y$} \\
      0  & \text{if $x \not\le y$}
   \end{cases}
\ee
(We think of $\zeta$ as a matrix whose rows and columns are indexed by $P$.)
The {\em M\"obius function}\/ $\mu$ on $P$ is the two-sided matrix inverse of
$\zeta$, i.e.\ it satisfies $\zeta\mu = \mu\zeta = I$.
It can be computed by the recursion
\be
   \mu(x,y)  \;=\;
   \begin{cases}
      1                                   & \text{if $x = y$} \\[1mm]
      -\sum\limits_{x \le z < y} \mu(x,z)  & \text{if $x < y$} \\[1mm]
      0                                   & \text{if $x \not\le y$}
   \end{cases}
\ee

Now let $S$ be a finite set,
and let $P$ be the lattice $\Pi(S)$ of partitions of $S$,
ordered by refinement.
It is well known \cite[p.~128]{Stanley_86} that the M\"obius function
of $\Pi(S)$ is given by
\be
   \mu(\sigma,\pi)  \;=\;
   \begin{cases}
       \prod\limits_{i=1}^k (-1)^{\lambda_i-1} (\lambda_i-1)!
             & \parbox{2.5in}{if $\sigma \le \pi = \{B_1,\ldots,B_k\}$
                        and $B_i$ is partitioned into $\lambda_i$ blocks
                        in $\sigma$}   \\[7mm]
       0     & \text{if $\sigma \not\le \pi$}
   \end{cases}
 \label{def.mu.partition}
\ee
We generalize this as follows:  let $q_1$ and $q_2$ be indeterminates,
and define
\be
   \mu_{q_1,q_2}(\sigma,\pi)  \;=\;
   \begin{cases}
       \prod\limits_{i=1}^k \prod\limits_{j=1}^{\lambda_i-1} (q_2 - jq_1)
             & \parbox{2.5in}{if $\sigma \le \pi = \{B_1,\ldots,B_k\}$
                        and $B_i$ is partitioned into $\lambda_i$ blocks
                        in $\sigma$}   \\[7mm]
       0     & \text{if $\sigma \not\le \pi$}
   \end{cases}
 \label{def.muq1q2.partition}
\ee
Let us observe for future reference that
\be
   \prod\limits_{j=1}^{m-1} (q_2 - jq_1)
   \;=\;
   q_1^m q_2^{-1} (q_2/q_1)^{\underline{m}}
   \;=\;
   q_1^{m-1} (q_2/q_1 - 1)^{\underline{m-1}}
   \;.
\ee
Then we have in particular
\begin{subeqnarray}
   \mu_{0,1}  & = &  \zeta \\[1mm]
   \mu_{1,0}  & = &  \mu \\[1mm]
   \mu_{q,q}  & = &  I \quad\hbox{for all $q$}
\end{subeqnarray}
We shall prove the following generalization of M\"obius inversion:

\begin{theorem}
  \label{thm.mobius}
Let $S$ be a finite set, and define
the matrices $\mu_{q_1,q_2}$ on $\Pi(S)$ by \reff{def.muq1q2.partition}.
Then
\be
   \mu_{q_1,q_2} \, \mu_{q_2,q_3}  \;=\; \mu_{q_1,q_3}
   \qquad\text{(matrix multiplication)}
 \label{eq.thm.mobius}
\ee
\end{theorem}

We interpret \reff{eq.thm.mobius} as an identity in the
polynomial ring $\Z[q_1,q_2,q_3]$,
but of course it also holds when $q_1,q_2,q_3$
are specialized to specific integer values
(or more generally to specific values in a commutative ring $R$).
In particular, M\"obius inversion corresponds to the special cases
$(q_1,q_2,q_3) = (0,1,0)$ and $(1,0,1)$.

The proof of Theorem~\ref{thm.mobius} will be based on the following lemma,
which we think is of some interest in its own right.
It is probably not new, but we have been unable to find any reference.

\begin{lemma}
   \label{lemma.partitions}
Let $m$ be a positive integer, and let $r$ and $s$ be indeterminates.  Then
\be
   \sum_{\omega \in \Pi_m}  r^{\underline{|\omega|}} \,
                     \Biggl( \prod_{B \in\omega} s^{\underline{|B|}} \Biggr)
   \;=\;\,
   (rs)^{\underline{m}}
 \label{eq.lemma.partitions}
\ee
as an identity in the polynomial ring $\Z[r,s]$.
\end{lemma}

Let us remark that if one divides \reff{eq.lemma.partitions}
by $s^m$ and takes $s \to\infty$
(or equivalently just extracts the coefficient of $s^m$),
one obtains the well-known formula
\reff{eq.stirlingidentity}/\reff{eq.stirling2};
while if one sets $r = q/\epsilon$,  $s = \epsilon q$
and takes $\epsilon\to 0$, then one obtains the formula
\be
   \sum_{\pi \in \Pi_m} \mu(\hat{0},\pi) \, q^{|\pi|}
   \;=\;
   q^{\underline{m}}
\ee
due to Rota \cite[section~9]{Rota_64}
\cite[pp.~128, 162, 187]{Stanley_86}
(here $\hat{0}$ denotes the partition in which every element is a singleton).

\proofof{Lemma~\ref{lemma.partitions}}
It suffices to prove \reff{eq.lemma.partitions} for positive integers $r,s$.
Then the right-hand side counts the proper colorings
of the complete graph $K_m$ with $rs$ colors.
But so does the left-hand side, if we use the color set
$[r] \times [s]$
and define $\omega$ to be the partition of $[m]$ in which two vertices
are placed in the same block if and only if
they receive a color with the same first index.
\qed

\bigskip\noindent
{\sc Second Proof of Lemma~\ref{lemma.partitions}\ }
(suggested independently by Christian Krattenthaler and Richard Stanley).
We use the isomorphism of weighted species
\be
   \hbox{Partitions}  \;=\; \hbox{Sets} \circ \hbox{NonemptySets} \;,
\ee
where a Set of cardinality $b$ is given a weight $r^{\underline{b}}$,
a NonemptySet of cardinality $k$ is given a weight $s^{\underline{k}}$,
and a partition $\omega$ is weighted as on the left-hand side
of \reff{eq.lemma.partitions}.\footnote{
   See \cite[pp.~44, 86]{Bergeron_98}
   for this isomorphism with a different weighting.
   For the reader unfamiliar with the theory of combinatorial species
   \cite{Bergeron_98}, it suffices to observe that a partition of a
   finite set $S$ is simply a set of nonempty subsets of $S$
   that disjointly cover $S$;
   one can then invoke \cite[Theorem~5.1.4]{Stanley_99}
   to complete the proof.
}
The corresponding exponential generating functions are
\begin{eqnarray}
   F(x) & = & \sum_{b=0}^\infty {x^b \over b!} \, r^{\underline{b}}
        \;=\; (1+x)^r  \\[1mm]
   G(x) & = & \sum_{k=1}^\infty {x^k \over k!} \, s^{\underline{k}}
        \;=\; (1+x)^s - 1
\end{eqnarray}
and hence
\be
   (F \circ G)(x) \;=\; (1+x)^{rs}
      \;=\; \sum_{m=0}^\infty {x^m \over m!} \, (rs)^{\underline{m}}
   \;,
\ee
which proves \reff{eq.lemma.partitions} for all $m \ge 0$.
\qed

\proofof{Theorem~\ref{thm.mobius}}
Let $\sigma \le \pi = \{B_1,\ldots,B_k\}$ where
$B_i$ is partitioned into $\lambda_i$ blocks in $\sigma$,
and let us compute
\be
   ( \mu_{q_1,q_2} \, \mu_{q_2,q_3} )(\sigma,\pi)
   \;=\;
   \sum_{\sigma \le \tau \le \pi} \mu_{q_1,q_2}(\sigma,\tau) \,
                                  \mu_{q_2,q_3}(\tau,\pi)   \;.
 \label{eq.proof.thm.mobius}
\ee
The partition $\tau$ is specified by saying, for each $i \in [k]$,
how the corresponding $\lambda_i$~blocks of $\sigma$ get grouped in $\tau$.
The sum \reff{eq.proof.thm.mobius} will then factorize over $i$,
with the $i$th factor given by
\begin{eqnarray}
   & &
   \sum_{\begin{scarray}
            \omega \in \Pi_{\lambda_i}  \\
            \omega = \{ \widehat{B}_1,\ldots, \widehat{B}_{\mu_i} \}
         \end{scarray}}
   \!\!\!
   \left( \prod\limits_{\alpha=1}^{\mu_i}
          \prod\limits_{j=1}^{|\widehat{B}_\alpha|-1} (q_2 -jq_1)
   \right)
   \prod\limits_{k=1}^{\mu_i-1} (q_3 -kq_2)
             \nonumber \\[2mm]
   & &
   \qquad\qquad =\;
   \!\!\!
   \sum_{\begin{scarray}
            \omega \in \Pi_{\lambda_i}  \\
            \omega = \{ \widehat{B}_1,\ldots, \widehat{B}_{\mu_i} \}
         \end{scarray}}
   \!\!\!
   \left( \prod\limits_{\alpha=1}^{\mu_i}
          q_1^{|\widehat{B}_\alpha|} q_2^{-1}
               (q_2/q_1)^{\underline{|\widehat{B}_\alpha|}}
   \right)
   q_2^{\mu_i} q_3^{-1} (q_3/q_2)^{\underline{\mu_i}}
             \nonumber \\[2mm]
   & &
   \qquad\qquad =\;
   q_1^{\lambda_i} q_3^{-1}
   \!\!\!
   \sum_{\begin{scarray}
            \omega \in \Pi_{\lambda_i}  \\
            \omega = \{ \widehat{B}_1,\ldots, \widehat{B}_{\mu_i} \}
         \end{scarray}}
   \!\!\!
   \left( \prod\limits_{\alpha=1}^{\mu_i}
               (q_2/q_1)^{\underline{|\widehat{B}_\alpha|}}
   \right)
   (q_3/q_2)^{\underline{\mu_i}}
             \nonumber \\[2mm]
   & &
   \qquad\qquad =\;
   q_1^{\lambda_i} q_3^{-1} (q_3/q_1)^{\underline{\lambda_i}}
\end{eqnarray}
by Lemma~\ref{lemma.partitions}.
\qed

\bigskip\noindent
{\sc Second Proof of Theorem~\ref{thm.mobius}\ }
(suggested by Richard Stanley).
We prove \reff{eq.thm.mobius} for $S=[n]$ simultaneously for all $n \ge 1$,
using the interpretation of convolution of multiplicative functions on
$\bm{\Pi} = (\Pi_1,\Pi_2,\ldots)$ in terms of composition of
exponential generating functions \cite[pp.~7--8]{Stanley_99}.
We first observe that $\mu_{q_1,q_2}$ is a multiplicative function
on $\bm{\Pi}$, with underlying numerical function
\be
   f_{q_1,q_2}(m)  \;=\;  \prod\limits_{j=1}^{m-1} (q_2 - jq_1)
                   \;=\;  q_1^m q_2^{-1} (q_2/q_1)^{\underline{m}}
\ee
and corresponding exponential generating function
\be
   F_{q_1,q_2}(x)  \;=\;  \sum_{m=1}^\infty f_{q_1,q_2}(m) \, {x^m \over m!}
                   \;=\;  q_2^{-1} [(1+ q_1 x)^{q_2/q_1} - 1]
                   \;.
\ee
It is easy to check that
\be
   F_{q_2,q_3}(F_{q_1,q_2}(x))  \;=\;  F_{q_1,q_3}(x)
   \;.
\ee
By \cite[Theorem~5.1.11]{Stanley_99},
the identities \reff{eq.thm.mobius} for all $n \ge 1$
are an immediate consequence.\footnote{
   Using \cite[Theorem~5.1.11]{Stanley_99} with the coefficient field $K$
   taken to be (say) $\R$, this argument proves \reff{eq.thm.mobius}
   for $q_1,q_2,q_3 \in \R \setminus \{0\}$,
   which is more than enough to imply it as a polynomial identity.
   Alternatively, we can use \cite[Theorem~5.1.11]{Stanley_99}
   with $K$ taken to be (say) the field $\R(q_1,q_2,q_3)$ of rational functions
   in the indeterminates $q_1,q_2,q_3$,
   which yields \reff{eq.thm.mobius} directly as a polynomial identity.
}
\qed

{\bf Remark.}
If we define the diagonal matrix $D_r$ by
\be
   D_r(\sigma,\pi)  \;=\;
   \begin{cases}
      r^{|\sigma|}   & \text{if $\sigma = \pi$} \\
      0              & \text{if $\sigma \neq \pi$}
   \end{cases}
\ee
then an easy calculation shows that
\be
   D_r \, \mu_{q_1,q_2} \, D_r^{-1}  \;=\;  \mu_{rq_1,rq_2}
   \;.
\ee
This is why many things involving $\mu_{q_1,q_2}$
depend only on the ratio $q_2/q_1$
[cf.\ \reff{eq.prop.ZG.partitions.q1q2}].
\qed

\bigskip

We can apply Theorem~\ref{thm.mobius} as follows:
Let $R$ be a commutative ring,
and fix a function $f \colon\, \Pi(S) \to R$
and an element $q_0 \in R$
(it could be an indeterminate if desired).
Now introduce a new indeterminate $q$,
and define $F_q \colon\, \Pi(S) \to R[q]$ by
\be
   F_q  \;=\; f \, \mu_{q_0,q}
\ee
or in more detail
\be
   F_q(\pi)  \;=\;
   \sum_{\sigma \in \Pi(S)}  f(\sigma) \, \mu_{q_0,q}(\sigma,\pi)
   \;.
\ee
By construction we have imposed the ``initial condition''
\be
   F_{q_0}  \;=\; f
\ee
(here $F_{q_0}$ means of course that the indeterminate $q$
 is replaced by the value $q_0 \in R$).
Most importantly, Theorem~\ref{thm.mobius} implies that
the $q$-dependence of $F_q$ is ``coherent'' in the sense that
\be
   F_{q_1} \, \mu_{q_1,q_2}  \;=\; F_{q_2}
   \;.
 \label{eq.Fq.coherent}
\ee

The formula \reff{eq.prop.ZG.partitions.q1q2.alt2}
is now the special case of \reff{eq.Fq.coherent}
in which we take $S=V$,
\be
   F_q(\pi)  \;=\;  \prod_{B\in \pi} \, \Zhat_{G[B]}(q,\bv)
   \;,
\ee
and we evaluate \reff{eq.Fq.coherent} at $\pi = \hat{1}$
(the partition in which all of $V$ belongs to a single block).
More generally, in place of $\Zhat_{G[B]}(q,\bv)$
we can use $\ahat_{\one_B}(q)$
where $\{ \ahat_\n(q) \}$ is {\em any}\/ family
arising as in Section~\ref{sec2};
we then obtain \reff{eq.partition2bis.q1hat.q2}
specialized to $\n = \one$.\footnote{
   Note also that, in these formulae,
   we use $\ahat_\n(q)$ [or $c_\n$] only for $\zero \le \n \le \one$.
}

\bigskip

The formulae in this appendix seem to be special for the partition lattice.
Still, it may not be totally absurd to ask:
Might there exist ``$q_1$--$q_2$ generalizations''
of M\"obius inversion for some other posets?
Maybe even for some natural {\em class}\/ of posets?

\section*{Acknowledgments}

We are grateful to Pierre Leroux for an extremely careful reading
of the manuscript and for many suggestions, both substantive and bibliographic.
In particular, it is thanks to him that our attention was drawn to
\cite[Section~3.1]{Bergeron_98}
and notably to \cite[p.~189, Exercise~3.1.22]{Bergeron_98},
from which we learned about the Abel-type extensions
\reff{eq.id1.t}--\reff{eq.id4.t}.

We finished the first version of this paper, by unhappy circumstance,
on the same day that we learned that
Pierre had been struck with a grave illness,
to which he succumbed two weeks later.
We dedicate this work to his memory.

Finally, we are grateful to an anonymous referee for pointing out
the connection of the formulae of Section~\ref{subsubsec.partitions}
with M\"obius inversion
--- which led us to discover the results reported in the Appendix ---
and to Christian Krattenthaler and Richard Stanley
for comments on these results.

This research was supported in part
by U.S.\ National Science Foundation grant PHY--0424082,
by U.K.\ Engineering and Physical Sciences Research Council grant GR/S26323/01,
and by a Leverhulme Visiting Professorship.



\end{document}